\documentclass[11pt]{article}
\usepackage[a4paper,margin=1in]{geometry}
\usepackage{amsmath,amssymb,amsthm,mathtools}
\usepackage[T1]{fontenc}
\usepackage[utf8]{inputenc}

\usepackage{esint}
\usepackage[colorlinks]{hyperref}
\usepackage[nameinlink,noabbrev]{cleveref}

\newtheorem{theorem}{Theorem}[section]
\newtheorem{lemma}[theorem]{Lemma}
\newtheorem{proposition}[theorem]{Proposition}

\theoremstyle{definition}
\newtheorem{definition}[theorem]{Definition}
\newtheorem{remark}[theorem]{Remark}

\title{Subdyadic time–frequency analysis: Gabor frames, modulation spaces, and Miyachi multipliers
}
\author{Vicente Vergara\footnote{Departamento de Matemática, Facultad de Ciencias Físicas y Matemáticas, Universidad de Concepción, Concepción, Chile. \texttt{vvergaraa@udec.cl}}}
\date{}

\begin{document}
	\maketitle
	
	\begin{abstract}
		We present a time–frequency framework adapted to dispersive phase functions via a subdyadic geometry in phase space. On top of this geometry we construct stable frequency-adaptive Gabor-type frames with quantitative control of overlap, almost orthogonality, and off-diagonal decay. Based on these frames we introduce modulation spaces consistent with the subdyadic scale and establish window and lattice independence, identifications in the Hilbertian case, duality, and natural inclusion relations. Within this setting we study high-frequency H\"ormander--Miyachi multipliers, relying on discrete block almost diagonalization and direct localization estimates for the primal and canonical dual frames, and obtain boundedness on weighted modulation spaces. Finally, we give a subdyadic Gabor-frame characterization of H\"ormander's classical local wavefront set and recover the standard microlocality and ellipticity properties of order-zero pseudodifferential operators. Taken together, these results provide a unified analytical framework for time--frequency analysis, dispersive multiplier theory, and local microlocal analysis in the subdyadic geometry.
	\end{abstract}
	


\section{Introduction}

The Gabor transform (short-time Fourier transform, STFT) plays a central
role in time--frequency analysis, as it provides a representation that is
simultaneously localized in space and frequency and quantitatively reflects
the uncertainty principle. Its formalism leads naturally to stable frames,
modulation spaces, and localized representations of functions and operators
\cite{Gabor1946,FeichtingerGrochenig1989,Grochenig2001}. More generally,
Gabor analysis provides a useful bridge between phase-space geometry and
discrete coefficient descriptions. This motivates the search for adapted
frame systems when the relevant geometry is no longer the uniform
time--frequency geometry of the classical STFT.

Dispersive phenomena governed by radial phases of the form $|\xi|^{\alpha}$, with $\alpha>0$, exhibit an intrinsic anisotropy in phase space that is not well captured by the classical dyadic paradigm. In the frequency region $|\xi|\sim 2^{j}$, the relevant scale for operators with essentially linear phase in $\xi$ naturally leads to cells $(\Delta x,\Delta \xi)\sim(2^{-j},2^{j})$. For the radial phase $|\xi|^{\alpha}$, the natural first-order radial
frequency scale is $|\xi|^{1-\alpha}$, with reciprocal spatial scale
$|\xi|^{\alpha-1}$; this leads to the \emph{subdyadic resolution} used
throughout the paper. This redistribution of time–frequency resolution has materialized in refined subdyadic tools (square functions and maximal operators) that have proved effective in weighted harmonic analysis \cite{BeltranBennett2017} and that are naturally connected to the theory of Miyachi multipliers \cite{Miyachi1980Wave,Miyachi1981Singular} and to the classical body of Littlewood–Paley inequalities and extrapolation \cite{Stein1993Harmonic}.

More precisely, let $\alpha,\beta\in\mathbb{R}$. We consider multipliers $m:\mathbb{R}^{d}\setminus\{0\}\to\mathbb{C}$ supported in
\[
\{ \xi\in\mathbb{R}^{d} : |\xi|^{\alpha} \geq 1\}.
\]
We say that $m$ satisfies the \emph{Miyachi condition} if
\begin{equation}\label{Cond:Miya:0}
	|D^{\gamma}m(\xi)| \;\lesssim\; |\xi|^{-\beta + |\gamma|(\alpha-1)},\qquad \text{for all }\xi \text{ with }|\xi|^{\alpha}\geq 1,
\end{equation}
for all multi-indices $\gamma\in\mathbb{N}^{d}$ with $|\gamma|\leq N_{d}$, where
\[
N_{d} = \Big\lfloor \frac{d}{2}\Big\rfloor + 1.
\]
The implicit constant may depend on $d,\alpha,\beta$ and $N_{d}$, but not on $\xi$.

A more flexible version, suitable for the analysis of Beltrán and Bennett \cite{BeltranBennett2017}, consists in requiring only a local quadratic control of the derivatives. Let $\mathcal{B}$ be a family of Euclidean balls in $\mathbb{R}^{d}$. For a ball $B\subset\mathbb{R}^{d}$ we write $r(B)$ for its radius and $\mathrm{dist}(B,0)$ for the distance from the center of $B$ to the origin. Given a multiplier $m$ supported in $\{|\xi|^{\alpha}\geq 1\}$, we say that it satisfies a \emph{Hörmander--Miyachi} type condition if
\begin{equation}\label{Cond:H-Miya:0}
	\sup_{B} \; \mathrm{dist}(B,0)^{\,\beta + (1-\alpha)|\gamma|}
	\left(
	\frac{1}{|B|} \int_{B} |D^{\gamma}m(\xi)|^{2}\,d\xi
	\right)^{1/2}
	<\infty,
\end{equation}
for all multi-indices $\gamma$ with $|\gamma|\leq N_{d}$, where the supremum is taken over all \emph{$\alpha$-subdyadic} balls $B$, that is,
\[
\mathrm{dist}(B,0)^{\alpha} \geq 1
\qquad\text{and}\qquad
r(B) \sim \mathrm{dist}(B,0)^{1-\alpha},
\]
where $r(B)$ denotes the radius of $B$. These balls provide a natural covering of the region $\{|\xi|^{\alpha}\geq 1\}$ by blocks whose radial scale is precisely $|\xi|^{1-\alpha}$.

The condition \eqref{Cond:H-Miya:0} is strictly weaker than the pointwise condition \eqref{Cond:Miya:0} above, but it still captures the same oscillatory and growth regime of the symbols. This is due to the fact that it measures the regularity of $m$ in an averaged sense over each ball, and the weight $\mathrm{dist}(B,0)^{\,\beta + (1-\alpha)|\gamma|}$ compensates the expected growth/decay of $D^{\gamma} m$.

In analogy with the $g$-functions developed by Stein \cite{Stein1970SingularIntegrals}, the authors of \cite{BeltranBennett2017} define subdyadic variants of the $g$-functions, constructed using convolution with a function $\varphi_t$. The corresponding \emph{subdyadic square function} is defined by
\[
g_{\alpha,\beta}(f)(x)
=
\left(
\int_{t^{\alpha}\leq 1}
\int_{|x-y|\leq t^{1-\alpha}}
|f\ast\varphi_{t}(y)|^{2}
\;t^{-(1-\alpha)d-2\beta}\,dy\,\frac{dt}{t}
\right)^{1/2}.
\]
For $\alpha=0$ and $\beta=0$ (with a suitable choice of $\varphi$) one recovers a variant of the classical Littlewood--Paley $g$-function. The parameter $\beta$ introduces a radial weight in the $t$ variable.

By analogy with the classical $g^{*}_{\lambda}$ function \cite{Stein1970SingularIntegrals}, one introduces
\[
g^{*}_{\alpha,\beta,\lambda}(f)(x)
=
\left(
\int_{t^{\alpha}\leq 1}
\int_{\mathbb{R}^{d}}
|f\ast\varphi_{t}(y)|^{2}
\left(1+\frac{|x-y|}{t^{1-\alpha}}\right)^{-d\lambda}
\;t^{-(1-\alpha)d-2\beta}\,dy\,\frac{dt}{t}
\right)^{1/2},
\]
with $\lambda>0$. By construction, $g^{*}_{\alpha,\beta,\lambda}$ is a pointwise majorant of $g_{\alpha,\beta}$, and for $\lambda$ sufficiently large it enjoys better analytic properties (for instance, in weighted estimates).

\subsubsection*{Two-sided Miyachi condition}\label{app:Miyachi-two}
Beltr\'an and Bennett \cite{BeltranBennett2017} also consider multipliers
defined on $\mathbb{R}^{d}\setminus\{0\}$ satisfying a two-sided version of
the Miyachi condition. In the present notation, for $\alpha>0$ this condition
takes the form
\begin{equation}\label{Cond:Miya:2-sided}
	|D^\gamma m(\xi)| \lesssim
	\begin{cases}
		|\xi|^{-\beta+(\alpha-1)|\gamma|},
		& |\xi|^\alpha\ge 1,\\[3pt]
		|\xi|^{-\beta-|\gamma|},
		& 0<|\xi|^\alpha\le 1,
	\end{cases}
\end{equation}
for all multi-indices $|\gamma|\le N_d$.
Thus the subdyadic behavior governed by $\alpha$ is imposed at high
frequency, whereas near the origin one uses the classical derivative scale.

The corresponding symbol class naturally includes the model symbols
\[
m_{\alpha,\beta}(\xi) = |\xi|^{-\beta}e^{\,i|\xi|^{\alpha}},
\]
which were studied in depth by Hirschman \cite{Hirschman1959}, Wainger \cite{Wainger1965}, Fefferman \cite{Fefferman1970}, Fefferman and Stein \cite{Fefferman1972}, Miyachi \cite{Miyachi1981Singular}, among others.

Beltr\'an and Bennett \cite{BeltranBennett2017} show that this two-sided
condition leads to weighted $L^2$ estimates through their subdyadic square
functions and maximal operators. In particular, their Corollary~11 gives the
following estimate.

\begin{theorem}[Corollary 11,\cite{BeltranBennett2017}]\label{teo:BB17:coro11:2-sided:miya}
	If $m:\mathbb{R}^{d}\setminus\{0\}\to\mathbb{C}$ satisfies \eqref{Cond:Miya:2-sided} for every $\gamma\in \mathbb{N}^d$ with $|\gamma|\leq N_d$, then
	\begin{equation}\label{Coro11:BB17:2-sided:miya}
		\int_{\mathbb{R}^{d}} |T_m f|^2 w \lesssim \int_{\mathbb{R}^{d}} |f|^2 M^2 \mathcal{M}_{\alpha, \beta} M^4 w,
	\end{equation}
	where
	\[
	\mathcal{M}_{\alpha, \beta}w(x) := \sup_{(y,r)\in \Lambda_{\alpha}(x)} \frac{1}{|B(y,r)|^{1-2\beta/d}} \int_{B(y,r)}w,
	\]
	and
	\[
	\Lambda_{\alpha}(x)
	=
	\{ (y,r)\in\mathbb{R}^{d}\times(0,\infty) : |x-y|\leq r^{\,1-\alpha}\}.
	\]
\end{theorem}
Here, $M$ denotes the classical Hardy–Littlewood maximal operator and $M^k$ stands for the $k$-fold composition of $M$ with itself. 

\medbreak

\subsection*{Relation to previous work and novelty}

This work proposes a \emph{subdyadic} time–frequency framework consistent with the phase $|\xi|^{\alpha}$ and tailored to microlocal applications. The starting point is an anisotropic metric on
$T^*\mathbb{R}^d_{\times}
=\mathbb R_x^d\times(\mathbb R_\xi^d\setminus\{0\})$
encoding the dispersive scale
\[
g^{\alpha}_{(x,\xi)}(dx,d\xi)\;\asymp\; |\xi|^{-2(\alpha-1)}|dx|^{2}+|\xi|^{-2(1-\alpha)}\,|d\xi|^{2},
\]
from which we define elementary blocks $Q_{\alpha}(x_0,\xi_0;1)$ of frequency size $\sim|\xi_0|^{1-\alpha}$ and spatial size $\sim|\xi_0|^{\alpha-1}$. The projection of $Q_{\alpha}$ onto the frequency space recovers the $\alpha$–subdyadic balls of \cite{BeltranBennett2017}, which allows us to inherit directly their partitions of unity and quadratic estimates. On top of this geometry we construct the frequency-dependent packet system
\[
\Omega_\alpha
=
\{(x_j^k,\xi_k):k\in\mathbb N,\ j\in\mathbb Z^d\},
\]
and prove that
$\{\phi_{x_j^k,\xi_k}\}_{(j,k)\in\Omega_\alpha}$
forms a stable $L^2$ frame. Since both the dilation of the window and the
spatial sampling step depend on the frequency index $k$, this is not a
classical single-window Gabor system. It belongs instead to the general
circle of frequency-adaptive or nonstationary Gabor constructions
\cite{DorflerMatusiak2014} and, from the frequency-covering viewpoint, to
the frame theory of decomposition spaces
\cite{BorupNielsen2007,Voigtlaender2022}.

We use the term \emph{subdyadic Gabor frame} to emphasize that the atoms
retain the local translation--modulation structure of Gabor analysis while
their resolution varies according to the subdyadic scale. This places the
construction naturally within the broader setting of frequency-adaptive
and nonstationary Gabor systems. The distinguishing feature of the present
system is the explicit choice
\[
\Delta\xi\asymp|\xi|^{1-\alpha},
\qquad
\Delta x\asymp|\xi|^{\alpha-1},
\]
dictated by the Beltr\'an--Bennett geometry
\cite{BeltranBennett2017}, together with its realization in the
phase-space metric $g^\alpha$ and the quantitative localization estimates
used later for dual frames, multiplier matrices, and microlocal analysis.

The explicit Fourier-side formula for the frame operator also gives direct
control of the canonical dual. In particular, the dual packets retain the
same subdyadic frequency supports as the primal packets and satisfy uniform
rescaled derivative bounds. This yields quantitative off-diagonal decay for
the dual and primal--dual Gram matrices by the same blockwise Fourier
localization argument used for the primal frame. The corresponding
Moore--Penrose identities are recorded in Lemma~\ref{lem:dual-local} and are
used later in the coefficient-space and multiplier analysis.

On the lattice $\Omega_\alpha$ we define coefficient spaces
$M^{p,q}_{\alpha,\beta}$ by means of $\ell^q(\ell^p)$ norms of the
subdyadic frame coefficients, with radial weight
$w_\beta(\xi_k)=(1+|\xi_k|)^\beta$. Their relation with the classical
theory depends on the value of $\alpha$. Set
\[
\theta:=1-\alpha.
\]
If $0<\alpha\le1$, then $0\le\theta<1$ and the projected frequency
blocks have diameter comparable to
$\langle\xi\rangle^\theta$. Hence they form a standard
$\theta$--modulation covering in the sense of
Gr{\"o}bner and of the subsequent Banach-frame theory
\cite{Groebner1992,BorupNielsen2006,Fornasier2007}. With the
$L^2$ normalization of the packets used here, the coefficient norm
corresponds to the usual $\theta$--modulation smoothness parameter
\[
s=\beta-\theta d\left(\frac12-\frac1p\right).
\]
Thus, in this parameter range, the spaces introduced below provide an
explicit subdyadic frame realization of the classical
$\theta$--modulation spaces.

When $\alpha>1$, the exponent $\theta=1-\alpha$ is negative and the
frequency cells shrink as $|\xi|\to\infty$. This regime lies outside
the standard $\theta$--modulation scale $0\le\theta\le1$, while the
resulting covering still fits naturally into the general framework of
structured decomposition spaces
\cite{FeichtingerGroebner1985,FeichtingerGroebner1987,
FeichtingerVoigtlaender2017,Voigtlaender2022}. The present construction
gives an explicit phase-space frame realization of this shrinking-cell
geometry, tied to the Beltr\'an--Bennett scale, and carries it into the
operator and microlocal estimates developed below.

Within this framework, we study Fourier multipliers satisfying the
high-frequency H\"ormander--Miyachi condition
\eqref{Cond:H-Miya:0}. Writing $\mu$ for the multiplier order, we prove
a block almost-diagonalization estimate for the corresponding operator
matrix with respect to the subdyadic frame. As a consequence,
\[
\|T_m f\|_{M^{p,q}_{\alpha,\beta}}
\lesssim
\|f\|_{M^{p,q}_{\alpha,\beta-\mu}},
\qquad
1\le p,q\le\infty.
\]
Thus the multiplier result holds on the full mixed-norm scale rather
than only in the Hilbertian case. Since
$M^{2,2}_{\alpha,\beta}=H^\beta$, the case $p=q=2$ reduces to the
corresponding Sobolev estimate. The subdyadic frame argument provides
the coefficient-level block estimate leading to the non-Hilbertian
mixed-norm bound, directly in the Beltr\'an--Bennett geometry and,
in particular, in the shrinking-cell regime $\alpha>1$.

Finally, we use the subdyadic frame to give an $\alpha$--dependent
Gabor-frame characterization of H\"ormander's classical local wavefront
set. For a spatially localized distribution, rapid decay of the
subdyadic frequency-block coefficient energies is shown to be equivalent
to classical conic Fourier regularity; consequently,
\[
WF_\alpha^{\mathrm G}(u)=WF(u),
\qquad \alpha>0.
\]
Thus the parameter $\alpha$ changes the analyzing resolution but not the
underlying set of local singularities. Standard pseudodifferential
microlocality and ellipticity then follow immediately from the classical
calculus.

This local construction is distinct from the global Gabor wavefront set
of Rodino and Wahlberg \cite{RodinoWahlberg2014}, which is defined by
conic decay in the full phase space and coincides with H\"ormander's
global wavefront set. Anisotropic global Gabor wavefront sets based on
power-type phase-space scaling were subsequently developed in
\cite{RodinoWahlberg2023,Wahlberg2024}. Related local wavefront
descriptions based on modulation spaces and discrete Gabor coefficients
can be found in \cite{JohanssonPilipovicTeofanovToft2012}.

Taken together, these results connect the Beltr\'an--Bennett subdyadic
geometry with a concrete phase-space frame and carry this realization
through function-space, multiplier, and microlocal analysis. The central
geometric feature is the simultaneous resolution
\[
\Delta\xi\asymp|\xi|^{1-\alpha},
\qquad
\Delta x\asymp|\xi|^{\alpha-1},
\]
implemented by an explicit Fourier-sampling frame construction. On this
frame, the multiplier analysis yields
\[
T_m:
M^{p,q}_{\alpha,\beta-\mu}
\longrightarrow
M^{p,q}_{\alpha,\beta},
\qquad
1\le p,q\le\infty,
\]
including the shrinking-frequency-cell regime $\alpha>1$. At the
microlocal level, the same frame provides an $\alpha$--dependent
characterization of H\"ormander's classical local wavefront set. In this
way, decomposition-space, nonstationary Gabor, and localized-frame
techniques are brought together within a phase-space geometry adapted
specifically to the subdyadic scale.

\medbreak

The paper is organized as follows. In Section
\ref{seccion:geoSubD:fase-espacio} we introduce the anisotropic metric
$g^{\alpha}$ on the punctured phase space, which encodes the natural scale
for dispersive phases of the form $|\xi|^{\alpha}$. We define the elementary
blocks $Q_{\alpha}(x,\xi;\rho)$ and establish their basic properties,
including their small-scale comparison with the balls of the distance
$d_{g^{\alpha}}$ for centers with $|\xi|\ge1$, the existence of coverings
with finite overlap (Lemma~\ref{lem:covering}), and uniform small-scale
volume and doubling estimates in the same frequency region
(Lemma~\ref{lem:volume}). Section \ref{seccion:RejillaGabor:STFTdispersiva} exploits this geometry to construct a lattice of points $\Omega_{\alpha}$ in phase space and to define the associated dispersive wave packets; we prove that the corresponding system $\mathcal{G}_{\alpha}(\phi)$ is an $L^2(\mathbb{R}^d)$ frame (Proposition~\ref{prop:frame}), with quantitative off-diagonal decay for the Gramian and the canonical dual (Proposition~\ref{prop:frame}). In Section \ref{seccion:EspaciosModulacion:Dispersivos}, on the subdyadic frame, we define a family of modulation spaces $M^{p,q}_{\alpha,\beta}$ (Definition~\ref{def:modulacion}) that measure the regularity and decay of a distribution in a way consistent with the dispersive dynamics, and we establish their basic properties: window and lattice independence (Theorem~\ref{thm:window-grid-indep}), identification of $M^{2,2}_{\alpha,\beta}$ with $H^\beta$
for every $\beta\in\mathbb R$
(Proposition~\ref{prop:L2}), density of Schwartz functions (Proposition~\ref{prop:densidad}), duality (Proposition~\ref{prop:dualidad}), and monotone inclusions (Proposition~\ref{prop:inclusiones}). 

Miyachi multipliers on the spaces $M^{p,q}_{\alpha,\beta}$ are studied
in Section \ref{seccion:Multiplicadores:Miyachi}. There we analyze
Fourier multipliers satisfying the high-frequency
H\"ormander--Miyachi condition, prove a block almost-diagonalization
estimate for their matrices with respect to the subdyadic frame
(Lemma~\ref{lem:matrix-Tm}), and deduce boundedness on the full scale
$M^{p,q}_{\alpha,\beta}$ for all $1\le p,q\le\infty$
(Theorem~\ref{thm:Miyachi-M22}).

Finally, Section~\ref{seccion:WFs:GaborSudDiadico} uses the subdyadic
frame coefficients to characterize H\"ormander's classical local
wavefront set. The resulting set $WF_\alpha^{\mathrm G}(u)$ is shown to
coincide with $WF(u)$ for every $\alpha>0$, and the standard
pseudodifferential microlocality and ellipticity properties are recovered
from the classical calculus.


\paragraph{Conventions and notation.}
We write $A\lesssim B$ if there exists a constant $C>0$
(independent of the relevant variables) such that $A\le C\,B$, and
$A\asymp B$ when both $A\lesssim B$ and $B\lesssim A$ hold.
We use $A\sim B$ informally to indicate comparability at the relevant
scale; when a quantitative two-sided estimate is needed in a proof, we
use $\asymp$. Unless otherwise stated, implicit constants are allowed to
depend only on fixed parameters of the problem (such as $d,\alpha$).

\section{Subdyadic geometry in phase space}\label{seccion:geoSubD:fase-espacio}

In this section we fix the \emph{reference geometry} in phase space $(x,\xi)\in \mathbb{R}^{d}_{x}\times\mathbb{R}^{d}_{\xi}$ that will be used to define subdyadic lattices, time–frequency blocks, and later on modulation spaces and wavefront sets.

\subsection{Natural scale for dispersive phases}

We consider dispersive phase functions of the form
\[
\varphi(\xi) = |\xi|^{\alpha}, \qquad \alpha>0,
\]
which give rise to evolution operators such as $e^{it(-\Delta)^{\alpha/2}}$ on $\mathbb{R}^{d}$ (see, for instance, \cite{KenigPonceVega1991, LinaresPonce2015}). The gradient has magnitude
\[
	|\nabla_{\xi}\varphi(\xi)|
	=
	\alpha|\xi|^{\alpha-1}.
\]
Hence, for increments $\Delta\xi$ with $|\Delta\xi|\ll|\xi|$, the
first-order variation satisfies
\[
	|\Delta\varphi(\xi)|
	\lesssim
	|\xi|^{\alpha-1}|\Delta\xi|.
\]
For radial increments this estimate is sharp up to constants, so the
radial frequency scale at which the phase changes by an amount of order
one is
\[
	|\Delta\xi|\asymp|\xi|^{1-\alpha}.
\]
This observation suggests that, in order to resolve accurately the oscillatory behavior of operators with phase $|\xi|^{\alpha}$, it is natural to decompose the frequency space into blocks whose characteristic length in the radial direction is comparable to $|\xi|^{1-\alpha}$. This is precisely the guiding principle behind the subdyadic decomposition introduced in \cite{BeltranBennett2017}, where frequency balls have radius tuned to $|\xi|^{1-\alpha}$ for frequencies of size $|\xi|$.

By Fourier duality, a frequency localization of size $|\Delta \xi|$ corresponds to a spatial localization of inversely proportional size, namely
\[
|\Delta x| \sim \frac{1}{|\Delta \xi|} \sim |\xi|^{\alpha-1}.
\]
In particular, the heuristic uncertainty product at this scale is given by
\[
|\Delta x|\,|\Delta \xi| \sim 1,
\]
but distributed in a non-uniform way as a function of $|\xi|$:
for $\alpha>1$ the frequency cells shrink and the spatial cells expand
as $|\xi|$ increases, whereas for $0<\alpha<1$ the frequency cells grow
sublinearly and the spatial cells shrink; for $\alpha=1$ both scales
remain constant. This redistribution is what makes the subdyadic geometry
``aligned'' with the uncertainty principle for phases of the form
$|\xi|^{\alpha}$.

\subsection{Subdyadic metric on the cotangent bundle}

The previous discussion suggests introducing an anisotropic Riemannian metric
$g^{\alpha}$ on the punctured cotangent bundle
\[
T^*\mathbb R^d_{\times}
:=\mathbb R^d_x\times(\mathbb R^d_\xi\setminus\{0\})
\]
of the form
\[
g^{\alpha}_{(x,\xi)}(dx,d\xi)
=|\xi|^{-2(\alpha-1)}|dx|^2+|\xi|^{-2(1-\alpha)}|d\xi|^2.
\]
Here $|dx|$ and $|d\xi|$ denote Euclidean norms on $\mathbb R^d$. This metric
assigns to a variation $d\xi$ of size $|d\xi|\sim|\xi|^{1-\alpha}$ a cost
comparable to that of a variation $dx$ of size $|\xi|^{\alpha-1}$, so that
its small balls in the high-frequency region reflect the dispersive scaling
discussed above. See \cite[Section 2.2]{Lerner2010} for background on
H\"ormander-type metrics.

\medskip
\noindent\textbf{Associated distance and balls.}
We denote by $d_{g^\alpha}$ the geodesic distance induced by $g^\alpha$ on
$T^*\mathbb R^d_{\times}$:
\[
d_{g^\alpha}\big((x,\xi),(y,\eta)\big)
:=\inf_\gamma\int_0^1
\Big(|\xi(t)|^{-2(\alpha-1)}|x'(t)|^2
+|\xi(t)|^{-2(1-\alpha)}|\xi'(t)|^2\Big)^{1/2}\,dt,
\]
where the infimum is taken over absolutely continuous curves
$\gamma(t)=(x(t),\xi(t))$ in $T^*\mathbb R^d_{\times}$ with prescribed
endpoints.  When $d=1$, the punctured frequency space
$\mathbb R\setminus\{0\}$ has two connected components. Accordingly,
$d_{g^\alpha}$ is understood as an extended geodesic distance on
$T^*\mathbb R_{\times}$, taking the value $+\infty$ between points whose
frequency components have opposite signs, and as the ordinary geodesic
distance on each connected component. All local comparisons below take
place within a single component.

For $z\in T^*\mathbb R^d_{\times}$ and $\rho>0$, let
\[
B_{g^\alpha}(z,\rho)
:=\{z'\in T^*\mathbb R^d_{\times}:d_{g^\alpha}(z,z')\le\rho\}.
\]

Later, in Lemma~\ref{lem:volume}, we show a uniform small-scale doubling
estimate for balls centered in the region $|\xi|\ge1$. For background on
doubling metric measure spaces, see \cite{Heinonen2001}.

\subsection{Subdyadic blocks in phase space}

For $w=(x,\xi), z=(y,\eta)\in T^*\mathbb R^d_{\times}$ we set
\[
r(w,z):=1+\min\{|\xi|,|\eta|\},
\]
and define the adapted separation function
\begin{equation}\label{def:quasi-d-alpha}
\mathbf d_\alpha(w,z)
:=r(w,z)^{1-\alpha}|x-y|+r(w,z)^{\alpha-1}|\xi-\eta|.
\end{equation}
We do not use $\mathbf d_\alpha$ as a global metric or quasi-metric; in
particular, no global triangle inequality is asserted. Its role is to record
the two subdyadic scales explicitly. For centers with $|\xi|\ge1$ and at
sufficiently small scale, the sets
\(\{z':\mathbf d_\alpha(z,z')\le C\rho\}\) are comparable, with uniform
constants, to the balls \(B_{g^\alpha}(z,\rho)\) and hence to the rectangular
blocks
\[
Q_\alpha(x,\xi;\rho)
=\{(y,\eta):|y-x|\le C\rho|\xi|^{\alpha-1},\
|\eta-\xi|\le C\rho|\xi|^{1-\alpha}\},
\]
as shown in Proposition~\ref{prop:distance} below.

This family of elementary blocks in phase space $(x,\xi)$ will serve as “building blocks’’ for the subsequent subdyadic Gabor decompositions. Note that these blocks are essentially the balls of $g^{\alpha}$, with the difference that we make the two scales explicit:
\begin{itemize}
	\item spatial scale: $|x-x_{0}|\sim \rho |\xi_0|^{\alpha-1}$;
	\item frequency scale: $|\xi-\xi_{0}|\sim \rho\,|\xi_{0}|^{1-\alpha}$.
\end{itemize}

When we fix $\rho\sim 1$, the blocks $Q_{\alpha}(x_{0},\xi_{0};1)$ have spatial size proportional to $|\xi_0|^{\alpha-1}$ and frequency size proportional to $|\xi_{0}|^{1-\alpha}$. By varying $\rho$, we can refine or enlarge the blocks in a way that is uniform with respect to the metric $g^{\alpha}$.

From the point of view of pure frequency decomposition (without the $x$ variable), the projection of $Q_{\alpha}(x_{0},\xi_{0};1)$ onto the $\xi$-variable is exactly a subdyadic ball centered at $\xi_{0}$ of radius comparable to $|\xi_{0}|^{1-\alpha}$. This connects directly with the subdyadic partitions of unity used in the analysis of Miyachi multipliers \cite{Miyachi1981Singular}.

\subsection*{Comparison with a dyadic reference}
In the standard dyadic decomposition by annuli
$\{2^{j-1}\le |\xi|\le 2^{j+1}\}$, a phase-space cell adapted to a
frequency block of size $2^j$ has the reciprocal scales
\[
	\Delta x \sim \rho\,2^{-j},
	\qquad
	\Delta \xi \sim \rho\,2^{j}.
\]
Accordingly, we may use the dyadic reference block
\[
	Q_{\mathrm{dy}}(x_0,\xi_0;\rho)
	:=
	\Big\{(x,\xi):\
	|x-x_0|\le C\rho\,2^{-j},\
	|\xi-\xi_0|\le C\rho\,2^{j}\Big\},
	\qquad |\xi_0|\sim 2^{j}.
\]
In contrast, for phases of the form $|\xi|^{\alpha}$ (with $\rho=1$), the natural scale of variation is no longer $2^{j}$ but $|\xi|^{1-\alpha}$, which in the annulus $|\xi|\sim 2^{j}$ translates into
\[
\Delta \xi \sim 2^{j(1-\alpha)}.
\]
This suggests that, at least for accurately describing the oscillatory behavior and almost orthogonality of dispersive operators, it is more efficient to “tilt’’ the decomposition towards blocks with this subdyadic scaling, rather than using isotropic blocks of size $2^{j}$ in frequency. The metric $g^{\alpha}$ provides a natural framework to formalize this idea \cite{BeltranBennett2017}. We will also make use of Littlewood–Paley type inequalities \cite{Stein1993Harmonic,RubioDeFrancia1985} and subdyadic partitions of unity \cite{BeltranBennett2017}.


\begin{proposition}[Local comparison with the $g^\alpha$ distance]\label{prop:distance}
	Let $\alpha>0$ and $z=(x,\xi)\in T^*\mathbb R^d_{\times}$ with
	$|\xi|\ge1$. There exist $c,C_1,C_2>0$, depending only on
	$d,\alpha$, such that, for every $\rho\in(0,c)$,
	\[
	B_{g^\alpha}(z,\rho)\;\subset\;
	\{z'\in T^*\mathbb R^d_{\times}:
	\mathbf d_\alpha(z,z')\le C_1\rho\}
	\;\subset\;
	B_{g^\alpha}(z,C_2\rho).
	\]
	In particular, for $\rho\le c$ the $g^\alpha$-balls are comparable,
	up to constants depending only on $d,\alpha$, to rectangular blocks
	of the form
	\[
	Q_\alpha(x,\xi;\rho):=\big\{(y,\eta):\,
	|y-x|\le C_3\rho|\xi|^{\alpha-1},\
	|\eta-\xi|\le C_3\rho|\xi|^{1-\alpha}\big\},
	\]
	for a fixed constant $C_3>0$ depending only on $d,\alpha$.
\end{proposition}

\begin{proof}
	Fix $z=(x,\xi)$ with $|\xi|\ge1$, and let $\rho>0$ be sufficiently
	small, depending only on $d$ and $\alpha$. All curves below are taken in
	$T^*\mathbb R^d_{\times}$.

	\medskip
	\noindent\emph{Step 1: Frequency comparability along short
	$g^\alpha$-paths.}
	Let $\gamma:[0,1]\to T^*\mathbb R^d_{\times}$,
	$\gamma(t)=(x_t,\xi_t)$, be an absolutely continuous path joining
	$z=(x,\xi)$ to $z'=(y,\eta)$, and suppose that
	\[
		L(\gamma)
		:=
		\int_0^1
		\sqrt{g^\alpha_{\gamma(t)}
		(\gamma'(t),\gamma'(t))}\,dt
		\le2\rho.
	\]
	From the definition of $g^\alpha$,
	\[
		\sqrt{g^\alpha_{\gamma(t)}
		(\gamma'(t),\gamma'(t))}
		\ge
		|\xi_t|^{\alpha-1}|\xi_t'|.
	\]
	Since
	\[
		\frac{d}{dt}|\xi_t|^\alpha
		=
		\alpha|\xi_t|^{\alpha-2}\xi_t\cdot\xi_t'
	\]
	for almost every $t$, we obtain
	\[
		\left|
		\frac{d}{dt}|\xi_t|^\alpha
		\right|
		\le
		\alpha|\xi_t|^{\alpha-1}|\xi_t'|.
	\]
	Hence
	\[
		\big||\xi_t|^\alpha-|\xi|^\alpha\big|
		\le
		\alpha L(\gamma)
		\le2\alpha\rho,
		\qquad 0\le t\le1.
	\]
	Choose $\rho_0>0$ so small that
	$2\alpha\rho_0\le1/2$. Since $|\xi|^\alpha\ge1$, for
	$\rho\le\rho_0$ this gives
	\[
		\frac12|\xi|^\alpha
		\le
		|\xi_t|^\alpha
		\le
		\frac32|\xi|^\alpha,
	\]
	and therefore
	\begin{equation}\label{eq:xi-comparable2}
		|\xi_t|\asymp|\xi|,
		\qquad 0\le t\le1,
	\end{equation}
	with constants depending only on $\alpha$. In particular,
	\[
		|\eta|\asymp|\xi|,
		\qquad
		r(z,z')
		=
		1+\min\{|\xi|,|\eta|\}
		\asymp|\xi|.
	\]

	\medskip
	\noindent\emph{Step 2: From the geodesic distance to the adapted
	separation function.}
	Assume that
	$d_{g^\alpha}(z,z')<\rho\le\rho_0$. Choose a path as above with
	$L(\gamma)\le2\rho$. Since
	\[
		\sqrt{a^2+b^2}
		\ge 2^{-1/2}(a+b),
		\qquad a,b\ge0,
	\]
	we have
	\[
		L(\gamma)
		\gtrsim
		\int_0^1
		\left(
		|\xi_t|^{1-\alpha}|x_t'|
		+
		|\xi_t|^{\alpha-1}|\xi_t'|
		\right)\,dt.
	\]
	Using \eqref{eq:xi-comparable2},
	\[
		L(\gamma)
		\gtrsim
		|\xi|^{1-\alpha}|x-y|
		+
		|\xi|^{\alpha-1}|\xi-\eta|.
	\]
	Since $r(z,z')\asymp|\xi|$, it follows that
	\[
		\mathbf d_\alpha(z,z')
		\lesssim
		L(\gamma)
		\lesssim\rho.
	\]
	Thus, after adjusting the constant,
	\[
		B_{g^\alpha}(z,\rho)
		\subset
		\{z':\mathbf d_\alpha(z,z')\le C\rho\}.
	\]

	\medskip
	\noindent\emph{Step 3: From the adapted separation function to the
	geodesic distance.}
	Conversely, suppose that
	\[
		\mathbf d_\alpha(z,z')\le\delta,
	\]
	where $\delta>0$ is sufficiently small. We first prove that
	\begin{equation}\label{eq:distance-endpoint-comparability}
		|\eta|\asymp|\xi|,
		\qquad
		r(z,z')\asymp|\xi|.
	\end{equation}

	If $|\eta|\ge|\xi|$, then
	$r(z,z')=1+|\xi|\asymp|\xi|$, and hence
	\[
		|\eta-\xi|
		\le
		\delta\,r(z,z')^{1-\alpha}
		\lesssim
		\delta|\xi|^{1-\alpha}.
	\]
	Consequently,
	\[
		|\eta|
		\le
		|\xi|+C\delta|\xi|^{1-\alpha}
		=
		|\xi|\bigl(1+C\delta|\xi|^{-\alpha}\bigr)
		\le
		(1+C\delta)|\xi|,
	\]
	so $|\eta|\asymp|\xi|$.

	If instead $|\eta|<|\xi|$, there is nothing to prove when
	$|\eta|\ge|\xi|/2$. Suppose therefore that
	$|\eta|<|\xi|/2$. Then
	\[
		\frac{|\xi|}{2}
		\le|\xi-\eta|
		\le
		\delta(1+|\eta|)^{1-\alpha}.
	\]
	If $\alpha\ge1$, then
	$(1+|\eta|)^{1-\alpha}\le1$, so
	$|\xi|/2\le\delta$, which is impossible for sufficiently small
	$\delta$ because $|\xi|\ge1$.
	If $0<\alpha<1$, then
	\[
		(1+|\eta|)^{1-\alpha}
		\lesssim|\xi|^{1-\alpha},
	\]
	and therefore
	\[
		\frac{|\xi|}{2}
		\lesssim
		\delta|\xi|^{1-\alpha},
	\]
	or equivalently
	\[
		|\xi|^\alpha\lesssim\delta,
	\]
	again impossible for sufficiently small $\delta$ since
	$|\xi|\ge1$. This proves
	\eqref{eq:distance-endpoint-comparability}.

	It follows from
	$\mathbf d_\alpha(z,z')\le\delta$ and
	\eqref{eq:distance-endpoint-comparability} that
	\begin{equation}\label{eq:distance-coordinate-bounds}
		|x-y|
		\lesssim
		\delta|\xi|^{\alpha-1},
		\qquad
		|\xi-\eta|
		\lesssim
		\delta|\xi|^{1-\alpha}.
	\end{equation}
	In particular,
	\[
		\frac{|\xi-\eta|}{|\xi|}
		\lesssim
		\delta|\xi|^{-\alpha}
		\le C\delta.
	\]
	Choosing $\delta$ still smaller if necessary, the linear path
	\[
		\gamma(t)
		=
		\bigl(x+t(y-x),\,\xi+t(\eta-\xi)\bigr),
		\qquad 0\le t\le1,
	\]
	stays in $T^*\mathbb R^d_{\times}$ and satisfies
	\[
		|\xi_t|\asymp|\xi|
		\qquad
		(0\le t\le1).
	\]
	Using \eqref{eq:distance-coordinate-bounds},
	\[
		\begin{aligned}
		L(\gamma)
		&=
		\int_0^1
		\sqrt{
		|\xi_t|^{2(1-\alpha)}|y-x|^2
		+
		|\xi_t|^{2(\alpha-1)}|\eta-\xi|^2
		}\,dt
		\\
		&\lesssim
		|\xi|^{1-\alpha}|x-y|
		+
		|\xi|^{\alpha-1}|\xi-\eta|
		\\
		&\lesssim
		\delta.
		\end{aligned}
	\]
	Thus
	\[
		d_{g^\alpha}(z,z')
		\le C\delta.
	\]
	Taking $\delta=C_0\rho$, with $C_0$ fixed, and decreasing the
	admissible upper bound for $\rho$ if necessary, gives
	\[
		\{z':\mathbf d_\alpha(z,z')\le C_0\rho\}
		\subset
		B_{g^\alpha}(z,C\rho).
	\]

	\medskip
	\noindent\emph{Step 4: Comparison with the rectangular blocks.}
	On either of the preceding small-scale sets we have
	$|\eta|\asymp|\xi|$, hence
	$r(z,z')\asymp|\xi|$. Therefore
	\[
		\mathbf d_\alpha(z,z')\lesssim\rho
	\]
	is equivalent, up to constants depending only on $d$ and $\alpha$, to
	\[
		|y-x|
		\lesssim
		\rho|\xi|^{\alpha-1},
		\qquad
		|\eta-\xi|
		\lesssim
		\rho|\xi|^{1-\alpha}.
	\]
	These are precisely the defining inequalities of
	$Q_\alpha(x,\xi;C\rho)$, up to a fixed change of constants.
	Combining this observation with Steps~2 and~3 yields
	\[
		B_{g^\alpha}(z,\rho)
		\subset
		Q_\alpha(x,\xi;C\rho)
		\subset
		B_{g^\alpha}(z,C'\rho),
	\]
	for all sufficiently small $\rho>0$. This proves the proposition.
\end{proof}

\begin{lemma}[Anisotropic volume and doubling property]\label{lem:volume}
	Let $\alpha>0$ and
	\[
	Q_\alpha(x_0,\xi_0;\rho)
	:=\big\{(y,\eta):\ |y-x_0|\le C\rho\,|\xi_0|^{\alpha-1},\ \ |\eta-\xi_0|\le C\rho\,|\xi_0|^{1-\alpha}\big\},
	\]
	with $|\xi_0|\ge 1$ and $\rho\in(0,\rho_0]$, where $\rho_0>0$ is sufficiently
	small and $C>0$ is a fixed constant. Then
	\begin{equation}\label{eq:vol-Qalpha}
		\operatorname{vol}\big(Q_\alpha(x_0,\xi_0;\rho)\big)\;\asymp\; \rho^{2d},
	\end{equation}
	with universal constants (independent of $(x_0,\xi_0)$ and $\rho$).
	
	In particular, geodesic balls centered in the region $|\xi|\ge1$ satisfy a
uniform small-scale doubling estimate: there exists
$C_{\mathrm{dbl}}=C_{\mathrm{dbl}}(d,\alpha)$ such that
	\begin{equation}\label{eq:doubling}
		\operatorname{vol}\big(B_{g^\alpha}(z,2\rho)\big)\;\le\;
		C_{\mathrm{dbl}}\ \operatorname{vol}\big(B_{g^\alpha}(z,\rho)\big),
		\qquad \rho\in(0,\rho_0],\ z\in\{(x,\xi):|\xi|\ge1\}.
	\end{equation}
\end{lemma}

\begin{proof}
	By definition, $Q_\alpha(x_0,\xi_0;\rho)$ is a rectangular block with side lengths
	\[
	\text{in the $x$ variable:}\quad \sim \rho\,|\xi_0|^{\alpha-1},\qquad
	\text{in the $\xi$ variable:}\quad \sim \rho\,|\xi_0|^{1-\alpha}.
	\]
	Since the reference measure is simply $dx\,d\xi$, we obtain
	\[
	\operatorname{vol}\big(Q_\alpha(x_0,\xi_0;\rho)\big)
	\asymp
	(\rho\,|\xi_0|^{\alpha-1})^{d}\,(\rho\,|\xi_0|^{1-\alpha})^{d}
	=\rho^{2d},
	\]
	which proves \eqref{eq:vol-Qalpha}. Note that the factor $|\xi_0|^{\alpha-1}$
	is compensated by $|\xi_0|^{1-\alpha}$, so that the volume is independent of $|\xi_0|$.
	
	For the doubling property, we use Proposition~\ref{prop:distance},
	which identifies the $d_{g^\alpha}$-balls with the blocks $Q_\alpha$. More precisely, there exist $c,C_1,C_2>0$ such that, for every $z=(x_0,\xi_0)$ with
	$|\xi_0|\ge1$ and every $\rho\in(0,\rho_0]$,
	\begin{equation}\label{eq:balls-Qalpha}
		B_{g^\alpha}(z,c\rho)\ \subset\ Q_\alpha(x_0,\xi_0;\rho)\ \subset\
		B_{g^\alpha}(z,C_1\rho),
	\end{equation}
	and similarly
	\[
	B_{g^\alpha}(z,\rho)\ \subset\ Q_\alpha(x_0,\xi_0;C_2\rho)\ \subset\
	B_{g^\alpha}(z,C_3\rho),
	\]
	after adjusting the constants.
	
	Then, using \eqref{eq:vol-Qalpha} and \eqref{eq:balls-Qalpha},
	\[
	\operatorname{vol}\big(B_{g^\alpha}(z,2\rho)\big)
	\ \le\
	\operatorname{vol}\big(Q_\alpha(x_0,\xi_0;C_4\rho)\big)
	\ \asymp\ (C_4\rho)^{2d}
	\ \lesssim\ \rho^{2d},
	\]
	while
	\[
	\operatorname{vol}\big(B_{g^\alpha}(z,\rho)\big)
	\ \gtrsim\
	\operatorname{vol}\big(Q_\alpha(x_0,\xi_0;c_1\rho)\big)
	\ \asymp\ (c_1\rho)^{2d}
	\ \gtrsim\ \rho^{2d}.
	\]
	Comparing these estimates, we obtain
	\[
	\operatorname{vol}\big(B_{g^\alpha}(z,2\rho)\big)
	\ \le\ C_{\mathrm{dbl}}\,
	\operatorname{vol}\big(B_{g^\alpha}(z,\rho)\big),
	\]
	with $C_{\mathrm{dbl}}$ depending only on $d$ and $\alpha$. This proves
	\eqref{eq:doubling} and completes the proof.
\end{proof}

\begin{lemma}[Covering with finite overlap]\label{lem:covering}
	Let $\alpha>0$ and fix $\rho\in(0,\rho_0]$, with $\rho_0>0$ sufficiently small. There exists a lattice
	\[
	\Omega_\alpha :=\{(x_j^k,\xi_k)\}_{k\in\mathbb{N},\,j\in\mathbb{Z}^d}\subset \mathbb{R}^d\times(\mathbb{R}^d\setminus\{0\})
	\]
	such that, with $Q_\alpha$ defined by
	\[
	Q_\alpha(x,\xi;\rho):=\big\{(y,\eta):\, |y-x|\le C\rho |\xi|^{\alpha-1},\ \ |\eta-\xi|\le C\rho\,|\xi|^{1-\alpha}\big\},
	\]
	the following properties hold:
	\begin{enumerate}
		\item[(1)] \emph{(Covering)}
\[
\{(x,\xi):|\xi|\ge1\}
\subset
\bigcup_{(j,k)}Q_\alpha(x_j^k,\xi_k;\rho).
\]
		\item[(2)] \emph{(Separation)} There exists $c_1\in(0,1)$ such that the blocks $Q_\alpha(x_j^k,\xi_k;c_1\rho)$ are pairwise disjoint.
		\item[(3)] \emph{(Finite overlap)} There exists $N_0=N_0(d,\alpha)$ such that every $(x,\xi)$ belongs to at most $N_0$ blocks $Q_\alpha(x_j^k,\xi_k;C\rho)$.
	\end{enumerate}
	Moreover, the lattice can be chosen by the following explicit procedure. For
\[
A_n:=\{\xi\in\mathbb{R}^d:2^n\le |\xi|<2^{n+1}\},
\qquad h_n:=\rho\,2^{n(1-\alpha)},\qquad n\in\mathbb{N}_0,
\]
take an $O(h_n)$--dense cubic mesh of candidate points in each $A_n$. Process these
candidates annulus by annulus and retain a candidate only when the Euclidean ball
centered there with radius $a\rho|\xi|^{1-\alpha}$ is disjoint from all such balls
previously retained, where $a>0$ is a fixed sufficiently small constant. Enumerate
the retained frequency centers as $\{\xi_k\}_{k\in\mathbb{N}}$, and let $n(k)$ be
the annulus containing $\xi_k$. For each $k$ set
\[
x_j^k:=b\rho|\xi_k|^{\alpha-1}j,\qquad j\in\mathbb{Z}^d,
\]
with a fixed $b>0$. The constants $a,b$ may be chosen depending only on $d$ and
$\alpha$, and the resulting covering, separation and overlap constants are independent
of $k$, $n(k)$ and $\rho$.
\end{lemma}

\begin{proof}
    We give the construction directly; it is a concrete realization of the maximal
    separated-net argument behind the Vitali/Besicovitch proof. For
    $n\in\mathbb{N}_0$ put
    \[
    A_n:=\{\xi:\ 2^n\le |\xi|<2^{n+1}\},
    \qquad h_n:=\rho\,2^{n(1-\alpha)}.
    \]
    Partition $\mathbb{R}^d$ into half-open cubes of side length $h_n$. From every cube
    meeting $A_n$ choose one point of the intersection. The resulting candidate set
    $\mathcal C_n\subset A_n$ is $C_dh_n$--dense in $A_n$, with $C_d$ depending only on
    $d$.

    \smallskip
    \noindent\emph{Step 1: Selection of the frequency centers.}
    Order the points of $\bigcup_{n\ge0}\mathcal C_n$ first by $n$ and then in any fixed
    order inside each annulus. Starting with the first point, retain a candidate $\xi$ if
    the Euclidean ball
    \[
    B_{\xi}^{\mathrm{in}}
    :=\big\{\eta:\ |\eta-\xi|<a\rho|\xi|^{1-\alpha}\big\}
    \]
    is disjoint from all inner balls corresponding to previously retained candidates.
    Here $a>0$ is fixed sufficiently small. Denote the retained points by
    $\{\xi_k\}_{k\in\mathbb N}$ and write $n(k)$ for the integer such that
    $\xi_k\in A_{n(k)}$.

    Since $\rho\le\rho_0$ and $\rho_0$ is small, a ball of radius
    $O(\rho|\xi|^{1-\alpha})$ centered at $\xi$ meets only the same or a neighboring
    dyadic annulus. Consequently, throughout such a ball the quantities
    $|\eta|^{1-\alpha}$ and $|\xi|^{1-\alpha}$ are comparable, with constants depending
    only on $\alpha$. The maximality of the preceding selection, together with the
    $C_dh_n$--density of the candidates, therefore gives a constant $A>1$ such that
    \[
    \{\eta:\ |\eta|\ge1\}
    \subset
    \bigcup_k
    \big\{\eta:\ |\eta-\xi_k|\le A\rho|\xi_k|^{1-\alpha}\big\}.
    \]
    The inner balls $B_{\xi_k}^{\mathrm{in}}$ are pairwise disjoint by construction.
    Moreover, the enlarged frequency balls above have uniformly bounded overlap. Indeed,
    if a fixed $\eta$ belongs to one of them, then $|\xi_k|\asymp|\eta|$; all
    corresponding pairwise disjoint inner balls are contained in a Euclidean ball centered
    at $\eta$ of radius $O(\rho|\eta|^{1-\alpha})$ and have radii comparable to
    $\rho|\eta|^{1-\alpha}$. A Euclidean volume comparison gives a bound depending only
    on $d$ and $\alpha$.

    \smallskip
    \noindent\emph{Step 2: Spatial grids.}
    For each retained frequency center define
    \[
    s_k:=b\rho|\xi_k|^{\alpha-1},
    \qquad
    x_j^k:=s_kj,\qquad j\in\mathbb Z^d,
    \]
    where $b>0$ is fixed. Every $x\in\mathbb R^d$ has a lattice point $x_j^k$ satisfying
    \[
    |x-x_j^k|\le \frac{\sqrt d}{2}\,s_k
    \lesssim \rho|\xi_k|^{\alpha-1}.
    \]
    Thus, after fixing the constant $C$ in the definition of $Q_\alpha$ sufficiently
    large, the frequency covering above and the preceding estimate imply
    \[
\{(x,\xi):|\xi|\ge1\}
\subset
\bigcup_{j,k}Q_\alpha(x_j^k,\xi_k;\rho).
\]

    Choose now $c_1>0$ sufficiently small. If $k=k'$ and $j\ne j'$, then
    $|x_j^k-x_{j'}^k|\ge s_k$, so the spatial projections of
    $Q_\alpha(x_j^k,\xi_k;c_1\rho)$ and
    $Q_\alpha(x_{j'}^k,\xi_k;c_1\rho)$ are disjoint. If $k\ne k'$, the frequency
    projections are disjoint once $c_1$ is chosen small compared with $a$, because the
    inner frequency balls selected in Step~1 are pairwise disjoint. This proves the
    separation property.

    Finally, fix $(x,\xi)$. By Step~1 there are only $O_{d,\alpha}(1)$ indices $k$ for
    which $|\xi-\xi_k|\lesssim\rho|\xi_k|^{1-\alpha}$. For each such $k$, the regular
    spatial lattice $s_k\mathbb Z^d$ contains only $O_d(1)$ points $x_j^k$ satisfying
    $|x-x_j^k|\lesssim\rho|\xi_k|^{\alpha-1}$. Hence $(x,\xi)$ belongs to at most
    $N_0=N_0(d,\alpha)$ enlarged blocks $Q_\alpha(x_j^k,\xi_k;C\rho)$. This proves the
    finite-overlap property and completes the construction.
\end{proof}

\medskip
\noindent\emph{Example ($d=1$ and $\alpha=2$).}
On the dyadic interval $[2^n,2^{n+1})$ the natural frequency mesh is
$h_n=\rho 2^{-n}$. One may therefore start with the candidates
\[
\xi_{n,m}^{+}=2^n+m h_n,
\qquad
\xi_{n,m}^{-}=-\bigl(2^n+m h_n\bigr),
\]
for the integers $m$ for which $2^n+m h_n<2^{n+1}$, and then apply the preceding
separation rule, discarding a candidate only when its inner frequency interval meets
one already retained. After enumerating the retained points as $\{\xi_k\}$, take
\[
x_j^k=b\rho|\xi_k|j,\qquad j\in\mathbb Z.
\]
Thus, near $|\xi|\sim2^n$,
\[
\Delta\xi\asymp\rho2^{-n},
\qquad
\Delta x\asymp\rho2^n,
\]
so that the phase-space cells have the reciprocal scales prescribed by $g^2$ and
$\Delta x\,\Delta\xi\asymp\rho^2$ uniformly in $n$.

\section{Subdyadic Gabor lattice and dispersive STFT}\label{seccion:RejillaGabor:STFTdispersiva}

\subsection{Construction of a subdyadic lattice in phase space}

Fix $\alpha>0$. Let $\rho_\ast\in(0,\rho_0]$ be an admissible scale in
Lemma~\ref{lem:covering}, to be chosen sufficiently small in
Proposition~\ref{prop:frame}. Once this choice is made, $\rho_\ast$ remains fixed
throughout the paper; we absorb it into the harmless geometric constants and, from
now on, write the corresponding blocks as $Q_\alpha(\cdot,\cdot;1)$. Apply the
explicit construction in Lemma~\ref{lem:covering}: the index $k$ enumerates the
retained frequency centers $\xi_k$, while $n(k)$ denotes the dyadic annulus
$2^{n(k)}\le|\xi_k|<2^{n(k)+1}$ containing $\xi_k$. The associated frequency balls
form an $\alpha$--subdyadic covering in the sense of \cite{BeltranBennett2017}, with
bounded overlap. Set
\[
r_k:=|\xi_k|^{\,1-\alpha},
\qquad
x_j^k=b\rho_\ast|\xi_k|^{\alpha-1}j,\qquad j\in\mathbb Z^d.
\]
Thus the frequency spacing near $\xi_k$ is $\asymp\rho_\ast r_k$ and the spatial
spacing is $\asymp\rho_\ast r_k^{-1}$. Equivalently, on the annulus
$|\xi|\sim2^n$ these scales are
\[
\Delta\xi\asymp2^{n(1-\alpha)},
\qquad
\Delta x\asymp2^{n(\alpha-1)},
\]
up to fixed constants. We therefore obtain the phase-space lattice
\[
\Omega_\alpha := \{z_{j,k}:=(x_j^k,\xi_k) : k\in\mathbb{N},\ j\in\mathbb{Z}^d\}
\subset \mathbb{R}^d\times(\mathbb{R}^d\setminus\{0\}),
\]
for which the blocks $Q_\alpha(x_j^k,\xi_k;1)$ cover
$\{(x,\xi):|\xi|^\alpha\ge1\}$, the smaller blocks
$Q_\alpha(x_j^k,\xi_k;c_1)$ are pairwise disjoint for some $c_1\in(0,1)$, and every
point belongs to at most $N_0$ enlarged blocks $Q_\alpha(x_j^k,\xi_k;C)$, with
$N_0$ independent of the annulus.

\subsection{Dispersive wave packets and adapted Gabor transform}

We fix a real-valued, even and $L^2$-normalized window
$\phi\in\mathcal{S}(\mathbb{R}^d)$ such that
$\widehat{\phi}\in C_c^\infty(\mathbb{R}^d)$. More precisely, we choose
constants $0<R_0<R_1<\infty$ and $m_0>0$, with
\[
R_0>3\,2^{\max\{\alpha-1,0\}},
\]
such that
\[
\operatorname{supp}\widehat{\phi}\subset B(0,R_1),
\qquad
\inf_{|\zeta|\le R_0}|\widehat{\phi}(\zeta)|\ge m_0.
\]
Such a window is obtained by taking a smooth even cutoff which is strictly
positive on $B(0,R_0)$ and then normalizing in $L^2$. For each point $(x_0,\xi_0)\in \mathbb{R}^d\times(\mathbb{R}^d\setminus\{0\})$ we define the \emph{dispersive wave packet}
\[
\phi_{x_0,\xi_0}(t)
=
|\xi_0|^{\frac{d}{2}(1-\alpha)}\,\phi\bigl((t - x_0)\,|\xi_0|^{1-\alpha}\bigr)\,e^{\,i\,t\cdot \xi_0},\qquad t\in\mathbb{R}^d.
\]
The normalization factor $|\xi_0|^{\frac{d}{2}(1-\alpha)}$ is chosen so that
$\|\phi_{x_0,\xi_0}\|_{L^2}$ is uniform in $(x_0,\xi_0)$, that is,
$\|\phi_{x_0,\xi_0}\|_{L^2} = \|\phi\|_{L^2}$. Note that $\phi_{x_0,\xi_0}$ is essentially concentrated in a spatial neighborhood of $x_0$ of radius $\sim|\xi_0|^{\alpha-1}$, in accordance with the spatial component of $Q_\alpha(x_0,\xi_0;1)$.

For brevity, we sometimes write $\phi_z$ to denote $\phi_{x,\xi}$, with $z=(x,\xi) \in \mathbb{R}^d\times(\mathbb{R}^d\setminus\{0\})$.

\begin{lemma}[Frequency localization of the packets]\label{lem:freq-local-packets}
	Suppose that $\widehat{\phi}$ is supported in $B(0,R_1)$. Then, for each
	$(x_0,\xi_0)$,
	\begin{equation}\label{eq:fourier-packet}
		\widehat{\phi_{x_0,\xi_0}}(\eta)
		=
		|\xi_0|^{-\frac d2(1-\alpha)}
		e^{-i x_0\cdot(\eta-\xi_0)}
		\widehat{\phi}\left(
		\frac{\eta-\xi_0}{|\xi_0|^{1-\alpha}}
		\right),
	\end{equation}
	and hence
	\[
	\operatorname{supp}\widehat{\phi_{x_0,\xi_0}}
	\subset
	\left\{
		\eta:\,
		|\eta-\xi_0|
		\le
		R_1|\xi_0|^{1-\alpha}
		\right\}.
	\]
	If, in addition,
	$|\widehat{\phi}(\zeta)|\ge m_0$ for $|\zeta|\le R_0$, then
	\[
	|\widehat{\phi_{x_0,\xi_0}}(\eta)|
	\ge
	m_0|\xi_0|^{-\frac d2(1-\alpha)}
	\]
	whenever
	\[
	|\eta-\xi_0|
	\le
	R_0|\xi_0|^{1-\alpha}.
	\]
\end{lemma}

\begin{proof}
	By the definition of the Fourier transform and the change of variables
	\[
	s=(t-x_0)|\xi_0|^{1-\alpha},
	\]
	we obtain
	\[
	\begin{aligned}
	\widehat{\phi_{x_0,\xi_0}}(\eta)
	&=
	|\xi_0|^{\frac d2(1-\alpha)}
	\int_{\mathbb{R}^d}
	\phi\big((t-x_0)|\xi_0|^{1-\alpha}\big)
	e^{it\cdot(\xi_0-\eta)}\,dt
	\\
	&=
	|\xi_0|^{-\frac d2(1-\alpha)}
	e^{-i x_0\cdot(\eta-\xi_0)}
	\widehat{\phi}\left(
	\frac{\eta-\xi_0}{|\xi_0|^{1-\alpha}}
	\right).
	\end{aligned}
	\]
	The support inclusion and the lower bound follow immediately from the
	corresponding properties of $\widehat{\phi}$.
\end{proof}

Thus, the packets $\phi_{x_0,\xi_0}$ are simultaneously
concentrated in a spatial neighborhood of $x_0$ of radius $\sim|\xi_0|^{\alpha-1}$ and in
a frequency neighborhood of $\xi_0$ of radius $\sim|\xi_0|^{1-\alpha}$, precisely at
the scale of the blocks $Q_\alpha(x_0,\xi_0;1)$.

\begin{remark}[Relation with nonstationary Gabor systems]
\label{rem:nonstationary-gabor}
	The terminology \emph{subdyadic Gabor system} used here does not mean
	a classical single-window Gabor system on a fixed rectangular lattice.
	Indeed, set
	\[
	r_k:=|\xi_k|^{1-\alpha},
	\qquad
	a_k:=b\rho_\ast r_k^{-1},
	\]
	and define on the Fourier side
	\[
	h_k(\eta)
	:=
	r_k^{-d/2}
	\widehat\phi\left(
	\frac{\eta-\xi_k}{r_k}
	\right).
	\]
	Then \eqref{eq:fourier-packet} gives
	\[
	\widehat{\phi_{x_j^k,\xi_k}}(\eta)
	=
	e^{i x_j^k\cdot\xi_k}\,
	e^{-i a_k j\cdot\eta}\,
	h_k(\eta).
	\]
	Thus, up to the irrelevant unimodular factor
	$e^{i x_j^k\cdot\xi_k}$, the Fourier-transformed system consists,
	for each frequency window $h_k$, of regular modulations with a
	$k$-dependent step $a_k$. Structurally this is a
	frequency-adaptive nonstationary Gabor system in the sense of the
	nonstationary Gabor-frame literature
	\cite{DorflerMatusiak2014}.

	The direct Fourier-sampling proof of
	Proposition~\ref{prop:frame} is correspondingly a
	multidimensional, frequency-side analogue of the painless
	nonstationary Gabor mechanism: compact frequency localization and
	sufficiently dense sampling reduce the frame inequalities to uniform
	upper and lower bounds for the local frequency density.

	More generally, the same packet architecture fits the structured
	frame constructions associated with decomposition coverings
	\cite{BorupNielsen2007,Voigtlaender2022}. We therefore do not claim
	the adaptive-frame principle itself as new. The feature specific to
	the present construction is that the window scale, frequency centers,
	and spatial sampling are all fixed explicitly by the
	Beltr\'an--Bennett subdyadic geometry.
\end{remark}

The \emph{dispersive Gabor transform} (or subdyadic STFT) of a function
$f\in L^2(\mathbb{R}^d)$ with respect to the window $\phi$ is defined by
\[
V_\phi f(x,\xi) :=\; \langle f, \phi_{x,\xi}\rangle
\;=\;
\int_{\mathbb{R}^d} f(t)\,\overline{\phi_{x,\xi}(t)}\,dt,
\qquad (x,\xi)\in\mathbb{R}^d\times(\mathbb{R}^d\setminus\{0\}).
\]
The function $V_\phi f(x,\xi)$ measures the correlation of $f$ with wave packets
localized in the region $Q_\alpha(x,\xi;1)$.

\begin{definition}[Separation and relative separation]\label{def:rel-sep}
	Let $(X,d)$ be a metric space. In our geometric application we take
$X=T^*\mathbb R^d_{\times}$ and $d=d_{g^\alpha}$.
	\begin{itemize}
		\item $\Lambda\subset X$ is \emph{$\delta$–separated} if $\inf_{\lambda\neq\mu\in\Lambda} d(\lambda,\mu)\ge\delta>0$.
		\item $\Lambda$ is \emph{relatively separated} if there exists $r_0>0$ such that
		\[
		N_{r_0}(\Lambda):=\sup_{z\in X}\#\big(\Lambda\cap B(z,r_0)\big)<\infty.
		\]
	\end{itemize}
	In \emph{doubling} spaces, this is equivalent (up to constants) to $\Lambda$ being a finite union of $r_0$–separated sets; see, for example, \cite[Chapter 11]{Grochenig2001} and also \cite[Section 2]{CorderoNicolaRodino2010}.
\end{definition}

\begin{lemma}[Relative separation of the lattice]\label{lem:rel-sep-Omega}
	Let $\alpha>0$ and $\rho\in(0,\rho_0]$. The lattice $\Omega_\alpha$ constructed in Lemma~\ref{lem:covering}
	is relatively separated in $(T^*\mathbb R^d_{\times},d_{g^\alpha})$.
More precisely, there exist	$r_0 (\asymp \rho)$ and $N_0=N_0(d,\alpha)$ such that
	\[
	\sup_{z\in T^*\mathbb R^d_{\times}} \#\big(\Omega_\alpha\cap B_{g^\alpha}(z,r_0)\big)\ \le\ N_0.
	\]
\end{lemma}

\begin{proof}
	By the \emph{separation} property in Lemma~\ref{lem:covering}, the blocks
	$Q_\alpha(\cdot;c_1\rho)$ are disjoint; by Proposition~\ref{prop:distance},
	$Q_\alpha(\cdot;c_1\rho)\supset B_{g^\alpha}(\cdot, c'\rho)$ for some $c'>0$.
	Set
\[
	r_0:=\frac{c'}{3}\rho.
\]
Then $r_0\asymp\rho$. If two distinct points
$\lambda,\mu\in\Omega_\alpha$ belonged to the same ball
$B_{g^\alpha}(z,r_0)$, the triangle inequality for the metric
$d_{g^\alpha}$ would give
\[
	d_{g^\alpha}(\lambda,\mu)
	\le
	d_{g^\alpha}(\lambda,z)+d_{g^\alpha}(z,\mu)
	<
	2r_0
	=
	\frac{2c'}{3}\rho
	<
	c'\rho,
\]
contradicting the separation established above. Therefore
\[
\sup_{z\in T^*\mathbb R^d_{\times}}
\#\bigl(\Omega_\alpha\cap B_{g^\alpha}(z,r_0)\bigr)
\le1.
\]
Thus $\Omega_\alpha$ is relatively separated, with $N_0=1$.
\end{proof}

\begin{proposition}[Subdyadic Gabor frame]\label{prop:frame}
	Let $\alpha>0$ and let $\phi$ be the window fixed above. Then there exists
	an admissible scale $\rho_\ast\in(0,\rho_0]$ such that the lattice
	$\Omega_\alpha$ constructed in Lemma~\ref{lem:covering} at scale $\rho_\ast$,
	with associated frequency centers $\{\xi_k\}$ and subdyadic balls $\{B_k\}$,
	has the following property. For the fixed window, satisfying
\[
\operatorname{supp}\widehat{\phi}\subset B(0,R_1),
\qquad
|\widehat{\phi}|\ge m_0 \quad\text{on }B(0,R_0),
\]
define the dispersive packets by
	\[
	\psi_{x,\xi}(t)= |\xi|^{\frac d2(1-\alpha)}\,\phi\bigl((t-x)\,|\xi|^{1-\alpha}\bigr)\,e^{i t\cdot \xi}.
	\]
	Then the system
	\[
	\mathcal{G}_\alpha(\phi)\;:=\;\{\psi_{x_j^k,\xi_k}:\ (x_j^k,\xi_k)\in\Omega_\alpha\}
	\]
	is a frame for $L^2(\mathbb{R}^d)$: there exist constants $0<A\le B<\infty$, depending only on $(d,\alpha)$, on $\phi$, and on the overlap constants of $\{B_k\}$ and $\Omega_\alpha$, such that
	\[
	A\,\|f\|_{L^2}^2\ \le\ \sum_{k}\sum_{j\in\mathbb{Z}^d}\big|\langle f,\psi_{x_j^k,\xi_k}\rangle\big|^2\ \le\ B\,\|f\|_{L^2}^2,\qquad f\in L^2(\mathbb{R}^d).
	\]
	Moreover, the canonical dual family $\{\widetilde{\psi}_{x_j^k,\xi_k}\}$ is intrinsically localized with respect to the same subdyadic phase-space geometry and provides the reconstruction formula
	\[
	f=\sum_k\sum_{j\in\mathbb{Z}^d} \langle f,\psi_{x_j^k,\xi_k}\rangle\,\widetilde{\psi}_{x_j^k,\xi_k},
	\quad\text{with convergence in }L^2.
	\]
\end{proposition}

\begin{proof}
	We write $\Omega_\alpha=\{w=(x_j^k,\xi_k)\}$ and
$\psi_w:=\psi_{x_j^k,\xi_k}$. We denote by
$G_{wz}:=\langle\psi_w,\psi_z\rangle$ the Gram-matrix entries.
The Fourier-sampling argument below will establish directly that the
coefficient map is bounded and bounded below on $L^2$; only after these
estimates are obtained do we introduce the bounded analysis operator and
the associated frame operator.
	
	\medskip
	
	\noindent\emph{Step 1: Off--diagonal decay of the Gramian.}
For $w=(x,\xi)$ and $z=(y,\eta)$, set
\[
	a:=|\xi|^{1-\alpha},
	\qquad
	b:=|\eta|^{1-\alpha},
	\qquad
	r:=1+\min\{|\xi|,|\eta|\}.
\]
We claim that, for every $N>0$,
\begin{equation}\label{eq:gram-decay}
	|G_{wz}|
	\le
	C_N\bigl(1+\mathbf d_\alpha(w,z)\bigr)^{-N},
\end{equation}
with $C_N$ independent of $w,z$.

By Lemma~\ref{lem:freq-local-packets}, if the Fourier supports of
$\psi_w$ and $\psi_z$ are disjoint, then $G_{wz}=0$. Assume therefore
that they intersect. Then
\[
	|\xi-\eta|
	\le
	R_1(a+b).
\]
Since $\alpha>0$, this overlap condition implies
$|\xi|\asymp|\eta|$ uniformly: at large frequency this follows from
$s^{1-\alpha}=o(s)$ as $s\to\infty$, while the remaining bounded
frequency region is absorbed into the constants. Consequently,
\begin{equation}\label{eq:gram-comparable-scales}
	a\asymp b\asymp r^{1-\alpha},
	\qquad
	r^{\alpha-1}|\xi-\eta|\lesssim1.
\end{equation}

Using Plancherel and \eqref{eq:fourier-packet}, followed by the change
of variables $\zeta=\xi+a u$, we obtain, up to the fixed Fourier
normalization constant,
\[
	G_{wz}
	=
	\Big(\frac{a}{b}\Big)^{d/2}
	e^{i y\cdot(\xi-\eta)}
	\int_{\mathbb R^d}
	e^{-i a(x-y)\cdot u}\,q_{wz}(u)\,du,
\]
where
\[
	q_{wz}(u)
	:=
	\widehat\phi(u)\,
	\overline{\widehat\phi\!\left(
	\frac{\xi-\eta}{b}+\frac{a}{b}u
	\right)}.
\]
By \eqref{eq:gram-comparable-scales}, the ratios $a/b$ and $b/a$ are
uniformly bounded, and so is $(\xi-\eta)/b$. Since
$\widehat\phi\in C_c^\infty(\mathbb R^d)$, for every integer $M\ge0$,
\[
	\big\|(1-\Delta_u)^M q_{wz}\big\|_{L^1_u}
	\le C_M,
\]
uniformly in $w,z$. Integrating by parts in $u$ therefore gives
\[
	|G_{wz}|
	\le
	C_M\bigl(1+a^2|x-y|^2\bigr)^{-M}.
\]
Since $a\asymp r^{1-\alpha}$, for every $N>0$ we may choose $M$
large enough to obtain
\[
	|G_{wz}|
	\le
	C_N\bigl(1+r^{1-\alpha}|x-y|\bigr)^{-N}.
\]
On the nonzero part of the Gramian, the frequency term in
$\mathbf d_\alpha(w,z)$ is uniformly bounded by
\eqref{eq:gram-comparable-scales}. Hence, using
\[
	\mathbf d_\alpha(w,z)
	=
	r^{1-\alpha}|x-y|
	+
	r^{\alpha-1}|\xi-\eta|,
\]
the preceding estimate yields \eqref{eq:gram-decay}. Together with the
disjoint-support case, this proves the claimed off-diagonal decay.
	
\medskip
\noindent\emph{Step 2: Frame bounds and frame operator by Fourier sampling.}
Let
\[
r_k:=|\xi_k|^{1-\alpha},
\qquad
a_k:=b\rho_\ast r_k^{-1},
\qquad
x_j^k=a_kj.
\]
Let $A_{\mathrm{cov}}>0$ be the constant in the frequency covering obtained
in Lemma~\ref{lem:covering}, so that
\[
\{\eta:\ |\eta|\ge1\}
\subset
\bigcup_k B(\xi_k,A_{\mathrm{cov}}\rho_\ast r_k).
\]

We now choose the admissible scale $\rho_\ast$ so small that
\begin{equation}\label{eq:rho-frame-choice}
\rho_\ast\le\rho_0,
\qquad
A_{\mathrm{cov}}\rho_\ast<R_0,
\qquad
b\rho_\ast R_1<\pi.
\end{equation}
Such a choice is always possible. Moreover, it does not destroy the covering:
Lemma~\ref{lem:covering} holds for every $0<\rho\le\rho_0$, and decreasing
$\rho$ simply refines both the frequency mesh
$\rho|\xi_k|^{1-\alpha}$ and the spatial mesh
$\rho|\xi_k|^{\alpha-1}$ in the explicit construction of that lemma.

By Lemma~\ref{lem:freq-local-packets},
\[
\widehat{\psi_{x_j^k,\xi_k}}(\eta)
=
r_k^{-d/2}
e^{-ia_kj\cdot(\eta-\xi_k)}
\widehat{\phi}\left(\frac{\eta-\xi_k}{r_k}\right).
\]
Since $\operatorname{supp}\widehat{\phi}\subset B(0,R_1)$ and
$b\rho_\ast R_1<\pi$, we have
\[
R_1r_k<\frac{\pi}{a_k}.
\]
Hence the support of the $k$-th frequency window is contained in a
fundamental cube for the exponential system
$\{e^{ia_kj\cdot\eta}\}_{j\in\mathbb{Z}^d}$.
Parseval's identity for Fourier series on that cube gives
\[
\sum_{j\in\mathbb{Z}^d}
|\langle f,\psi_{x_j^k,\xi_k}\rangle|^2
=
(2\pi)^{-d}a_k^{-d}r_k^{-d}
\int_{\mathbb{R}^d}
|\widehat f(\eta)|^2
\left|
\widehat{\phi}\left(\frac{\eta-\xi_k}{r_k}\right)
\right|^2
\,d\eta.
\]
Since $a_k=b\rho_\ast r_k^{-1}$,
\[
a_k^{-d}r_k^{-d}=(b\rho_\ast)^{-d},
\]
independently of $k$. Summing over $k$ therefore yields
\begin{equation}\label{eq:painless-frame-identity}
\sum_{k,j}
|\langle f,\psi_{x_j^k,\xi_k}\rangle|^2
=
(2\pi)^{-d}(b\rho_\ast)^{-d}
\int_{\mathbb{R}^d}
|\widehat f(\eta)|^2P(\eta)\,d\eta,
\end{equation}
where
\[
P(\eta)
:=
\sum_k
\left|
\widehat{\phi}\left(\frac{\eta-\xi_k}{r_k}\right)
\right|^2.
\]

We claim that $P$ is bounded above and below by positive constants.
If $|\eta|\ge1$, the frequency covering provides an index $k$ such that
\[
|\eta-\xi_k|
\le
A_{\mathrm{cov}}\rho_\ast r_k
<
R_0r_k.
\]
Consequently,
\[
P(\eta)\ge m_0^2.
\]
If $|\eta|<1$, choose any center $\xi_{k_0}$ in the first dyadic annulus
$1\le|\xi_{k_0}|<2$. Since
\[
r_{k_0}=|\xi_{k_0}|^{1-\alpha},
\]
we have
\[
\frac{|\eta-\xi_{k_0}|}{r_{k_0}}
\le
3\,2^{\max\{\alpha-1,0\}}
<
R_0,
\]
and hence again $P(\eta)\ge m_0^2$.

For the upper bound, a summand in $P(\eta)$ can be nonzero only if
\[
|\eta-\xi_k|\le R_1r_k.
\]
The separation of the frequency centers from
Lemma~\ref{lem:covering}, together with a Euclidean packing argument on
the relevant dyadic annuli, implies
\[
\sup_{\eta\in\mathbb{R}^d}
\#\{k:\ |\eta-\xi_k|\le R_1r_k\}
\le N_1<\infty,
\]
where $N_1$ depends only on $d,\alpha,R_1$ and the fixed admissible
scale $\rho_\ast$. Thus
\[
m_0^2
\le
P(\eta)
\le
N_1\|\widehat{\phi}\|_{L^\infty}^2,
\qquad
\eta\in\mathbb{R}^d.
\]
Combining this estimate with
\eqref{eq:painless-frame-identity} and Plancherel's theorem gives
\[
(b\rho_\ast)^{-d}m_0^2\|f\|_{L^2}^2
\le
\sum_{k,j}
|\langle f,\psi_{x_j^k,\xi_k}\rangle|^2
\le
(b\rho_\ast)^{-d}
N_1\|\widehat{\phi}\|_{L^\infty}^2
\|f\|_{L^2}^2.
\]
Thus the coefficient map extends to a bounded operator
\[
C:L^2(\mathbb R^d)\longrightarrow\ell^2(\Omega_\alpha),
\qquad
Cf=(\langle f,\psi_w\rangle)_w,
\]
and the preceding two-sided estimate shows that
$\mathcal G_\alpha(\phi)$ is a frame for $L^2(\mathbb R^d)$.

Let $S=C^\ast C$ be the frame operator. Polarizing
\eqref{eq:painless-frame-identity}, for $f,g\in L^2(\mathbb R^d)$ we obtain
\[
\langle Sf,g\rangle
=
(2\pi)^{-d}(b\rho_\ast)^{-d}
\int_{\mathbb R^d}
\widehat f(\eta)\,\overline{\widehat g(\eta)}\,P(\eta)\,d\eta.
\]
With the Fourier convention used here, Plancherel's identity therefore gives
\begin{equation}\label{eq:frame-operator-symbol}
\widehat{Sf}(\eta)
=
(b\rho_\ast)^{-d}P(\eta)\widehat f(\eta).
\end{equation}
Since $P$ is bounded above and below by positive constants, $S$ is boundedly
invertible on $L^2(\mathbb R^d)$ and
\[
\widehat{S^{-1}f}(\eta)
=
(b\rho_\ast)^d P(\eta)^{-1}\widehat f(\eta).
\]
Hence the canonical dual packet associated with
$w=(x_j^k,\xi_k)$ satisfies
\begin{equation}\label{eq:dual-fourier-packet}
\widehat{\widetilde\psi_{x_j^k,\xi_k}}(\eta)
=
r_k^{-d/2}e^{-ix_j^k\cdot(\eta-\xi_k)}
A_k\!\left(\frac{\eta-\xi_k}{r_k}\right),
\end{equation}
where
\[
A_k(u)
:=
(b\rho_\ast)^d
P(\xi_k+r_k u)^{-1}\widehat\phi(u).
\]
In particular,
$\operatorname{supp}A_k\subset\operatorname{supp}\widehat\phi$.
Moreover, for every multi-index $\gamma$,
\begin{equation}\label{eq:dual-amplitude-uniform}
\sup_k\sup_{u\in\mathbb R^d}
|\partial_u^\gamma A_k(u)|
\le C_\gamma.
\end{equation}
Indeed, on $\operatorname{supp}\widehat\phi$ the point
$\xi_k+r_k u$ belongs to the $k$-th frequency support. Every term in
$P(\xi_k+r_k u)$ which is nonzero there comes from a frequency support
intersecting the $k$-th one. By \eqref{eq:gram-comparable-scales}, the
corresponding scales are uniformly comparable to $r_k$, while the packing
bound used above gives a uniformly bounded number of such terms. Therefore,
for every multi-index $\gamma$,
\[
\sup_k\sup_{u\in\operatorname{supp}\widehat\phi}
\left|
\partial_u^\gamma P(\xi_k+r_k u)
\right|
\le C_\gamma.
\]
Since $P\ge m_0^2$, the same bounds hold for
$P(\xi_k+r_k u)^{-1}$, and \eqref{eq:dual-amplitude-uniform} follows.

Formula \eqref{eq:dual-fourier-packet} shows that the canonical dual packets
have the same subdyadic Fourier supports as the primal packets and form a
uniformly localized family at the same scales. More precisely, repeating the
Fourier-side integration-by-parts argument of Step~1, now using the uniform
bounds \eqref{eq:dual-amplitude-uniform}, gives, for every $N>0$,
\begin{equation}\label{eq:dual-direct-localization}
|\langle\widetilde\psi_w,\widetilde\psi_z\rangle|
+
|\langle\widetilde\psi_w,\psi_z\rangle|
\le
C_N\bigl(1+\mathbf d_\alpha(w,z)\bigr)^{-N}.
\end{equation}
Indeed, if the corresponding Fourier supports are disjoint, both inner
products vanish; if they intersect, the scales are comparable and the
frequency component of $\mathbf d_\alpha(w,z)$ is uniformly bounded by
\eqref{eq:gram-comparable-scales}, while integration by parts yields
arbitrary decay in the scaled spatial separation.

Thus the canonical dual family is intrinsically localized with respect to the
same subdyadic phase-space geometry. Since it is the canonical dual of the
frame, the standard frame identity yields
\[
f
=
\sum_k\sum_{j\in\mathbb Z^d}
\langle f,\psi_{x_j^k,\xi_k}\rangle
\widetilde\psi_{x_j^k,\xi_k},
\]
with convergence in $L^2(\mathbb R^d)$.
\end{proof}

\section{Dispersive modulation spaces}
\label{seccion:EspaciosModulacion:Dispersivos}

In this section we introduce the dispersive modulation spaces $M^{p,q}_{\alpha,\beta}$ associated with the $\alpha$-subdyadic geometry and to the dispersive Gabor transform constructed in the previous section. The parameter $\alpha>0$ fixes the metric $g^\alpha$ in phase space, while $\beta\in\mathbb{R}$ controls the radial weight in frequency.

Recall that, for a fixed window $\phi\in\mathcal{S}(\mathbb{R}^d)$, the dispersive STFT of $f\in L^2(\mathbb{R}^d)$ was defined by
\[
V_\phi f(x,\xi)=\langle f,\phi_{x,\xi}\rangle,\qquad (x,\xi)\in\mathbb{R}^d\times(\mathbb{R}^d\setminus\{0\}),
\]
where the packets $\phi_{x,\xi}$ are localized in blocks $Q_\alpha(x,\xi;1)$ for the metric $g^\alpha$. We also fix a subdyadic Gabor lattice
\[
\Omega_\alpha=\{(x_j^k,\xi_k):\ k\in\mathbb{N},\ j\in\mathbb{Z}^d\},
\]
such that the blocks $Q_\alpha(x_j^k,\xi_k;1)$ cover $\{(x,\xi):|\xi|^\alpha\ge 1\}$ with finite overlap and the family
\[
\mathcal{G}_\alpha(\phi)=\{\phi_{x_j^k,\xi_k}:(x_j^k,\xi_k)\in\Omega_\alpha\}
\]
forms a frame for $L^2(\mathbb{R}^d)$ (Proposition~\ref{prop:frame}).

\paragraph{Low frequencies.}
Although the frequency centers of $\Omega_\alpha$ satisfy $|\xi_k|\ge1$,
no additional low-frequency atoms are needed. Indeed, the proof of
Proposition~\ref{prop:frame} shows that the windows associated with the first
dyadic annulus already yield
\[
P(\eta)\ge m_0^2,
\qquad |\eta|<1.
\]
Thus the same subdyadic frame covers the low-frequency part of $L^2$ as well,
without any finite modification of the system.

\subsection{Discrete definition via Gabor coefficients}

For each $f\in\mathcal{S}'(\mathbb{R}^d)$ we define the dispersive Gabor coefficients
\[
c_{j,k}(f)=\langle f,\phi_{x_j^k,\xi_k}\rangle,\qquad (j,k)\in\mathbb{Z}^d\times\mathbb{N}.
\]
We introduce the radial weight
\[
w_\beta(\xi)=(1+|\xi|)^\beta.
\]

\begin{definition}[Dispersive modulation spaces]\label{def:modulacion}
	Let $1\le p,q\le\infty$ and $\beta\in\mathbb{R}$. We define $M^{p,q}_{\alpha,\beta}$ as the set of $f\in\mathcal{S}'(\mathbb{R}^d)$ such that
	\[
	\|f\|_{M^{p,q}_{\alpha,\beta}}
	:=
	\left(
	\sum_k
	\Big(\sum_{j\in\mathbb{Z}^d}|c_{j,k}(f)|^p\Big)^{\!q/p}
	\,w_\beta(\xi_k)^q
	\right)^{\!1/q}
	<\infty,
	\]
	with the usual modifications when $p$ or $q$ is infinite (replacing sums by suprema).
\end{definition}

The anisotropy of $g^\alpha$ is reflected in the fact that the cells
$Q_\alpha(x_j^k,\xi_k;1)$ have radii in $x$ and $\xi$ that depend on
$|\xi_k|$ (of order $|\xi_k|^{\alpha-1}$ and $|\xi_k|^{1-\alpha}$,
respectively), although their volume is comparable to $1$ when $\rho=1$
(Lemma~\ref{lem:volume}).

\begin{proposition}[Relation with $\theta$--modulation spaces]
\label{prop:relation-alpha-mod}
	Assume $0<\alpha\le1$ and set
	\[
		\theta:=1-\alpha\in[0,1).
	\]
	For the frequency centers of Proposition~\ref{prop:frame}, put
	\[
		r_k:=|\xi_k|^{1-\alpha}
	\]
	and set
\[
	U_k:=B(\xi_k,R_1r_k),
\]
so that $U_k$ contains the Fourier support of the $k$-th packet family
with a uniform radius proportional to $r_k$.
	Then $(U_k)_k$ is, up to weak equivalence in the bounded-frequency
	region, an inhomogeneous $\theta$--modulation covering.

	For $1\le p,q\le\infty$ and $\beta\in\mathbb R$, set
	\[
		s:=\beta-\theta d\left(\frac12-\frac1p\right).
	\]
	Then
	\[
		M^{p,q}_{\alpha,\beta}
		=
		M^{s,\theta}_{p,q}(\mathbb R^d)
	\]
	with equivalence of norms.
\end{proposition}

\begin{proof}
	Set
	\[
		\Phi_k(\eta)
		:=
		\widehat\phi\!\left(\frac{\eta-\xi_k}{r_k}\right),
		\qquad
		P(\eta):=\sum_k|\Phi_k(\eta)|^2.
	\]
	By Proposition~\ref{prop:frame}, $P$ is bounded above and below by
	positive constants. Hence
	\[
		\chi_k(\eta):=\frac{|\Phi_k(\eta)|^2}{P(\eta)}
	\]
	defines a smooth partition of unity. Moreover,
	$\operatorname{supp}\chi_k\subset U_k$, and the rescaled derivatives of
	$\chi_k$ are bounded uniformly in $k$. Indeed, this is the same local
	overlap argument used in the proof of
	\eqref{eq:dual-amplitude-uniform}. Consequently
	\[
		\sup_k\big\|\mathcal F^{-1}\chi_k\big\|_{L^1}<\infty,
	\]
	so $(\chi_k)_k$ is a bounded admissible partition of unity for the
	covering $(U_k)_k$.

	Since $r_k\asymp\langle\xi_k\rangle^\theta$ and intersecting frequency
	sets have comparable centers and scales, the high-frequency part of
	$(U_k)_k$ has precisely the geometry of a $\theta$--modulation covering.
	Only finitely many sets are involved in the bounded-frequency region.
	Thus $(U_k)_k$ is weakly equivalent to a standard inhomogeneous
	$\theta$--modulation covering, and the associated decomposition norm
	\[
		\left\|
		\left(
		\langle\xi_k\rangle^s
		\big\|\mathcal F^{-1}(\chi_k\widehat f)\big\|_{L^p}
		\right)_k
		\right\|_{\ell^q}
	\]
	is an equivalent norm for $M^{s,\theta}_{p,q}(\mathbb R^d)$; see
	\cite{Groebner1992,BorupNielsen2006,Voigtlaender2022}.

	We now compare this decomposition norm directly with the coefficients of
	Definition~\ref{def:modulacion}. Put
	\[
		h_k
		:=
		\mathcal F^{-1}\!\left(
		\widehat f\,\overline{\Phi_k}
		\right).
	\]
	By Lemma~\ref{lem:freq-local-packets}, with
	$a_k=b\rho_\ast r_k^{-1}$,
	\[
		|c_{j,k}(f)|
		=
		r_k^{-d/2}|h_k(a_kj)|.
	\]
	Since $a_kr_k=b\rho_\ast$ is fixed and
	$b\rho_\ast R_1<\pi$, the usual sampling inequalities for functions
	whose Fourier transforms are supported in a fixed compact subset of a
	fundamental cube, after modulation and dilation to unit scale, give
	\begin{equation}\label{eq:theta-sampling-block}
		\big\|(c_{j,k}(f))_j\big\|_{\ell^p}
		\asymp
		r_k^{d(1/p-1/2)}\,\|h_k\|_{L^p},
		\qquad 1\le p\le\infty,
	\end{equation}
	with constants independent of $k$. The same assertion holds for
	bandlimited distributions: membership of the sampled sequence in
	$\ell^p$ is equivalent to the corresponding bandlimited distribution
	having an $L^p$ representative. One way to see this is to insert a smooth
	cutoff equal to one on the fixed rescaled Fourier support; Fourier-series
	reconstruction and the Schur test then give both inequalities in
	\eqref{eq:theta-sampling-block}.

	It remains to compare $h_k$ with the decomposition pieces. Since
	\[
		\chi_k
		=
		\frac{\Phi_k}{P}\,\overline{\Phi_k},
	\]
	the multiplier $\Phi_k/P$ has a uniformly bounded $L^1$ convolution
	kernel after rescaling. Hence
	\[
		\big\|\mathcal F^{-1}(\chi_k\widehat f)\big\|_{L^p}
		\lesssim
		\|h_k\|_{L^p}.
	\]
	Conversely, using $\sum_\ell\chi_\ell=1$,
	\[
		\widehat f\,\overline{\Phi_k}
		=
		\sum_{\ell:\,U_\ell\cap U_k\neq\varnothing}
		\overline{\Phi_k}\,\chi_\ell\widehat f.
	\]
	Only uniformly finitely many $\ell$ occur, and on interacting blocks the
	weights $\langle\xi_\ell\rangle^s$ and
	$\langle\xi_k\rangle^s$ are comparable. The multipliers
	$\overline{\Phi_k}$ on these blocks again have uniformly bounded
	$L^1$ convolution kernels. A finite-degree Schur estimate in the outer
	frequency index therefore yields
	\[
		\left\|
		\left(
		\langle\xi_k\rangle^s\|h_k\|_{L^p}
		\right)_k
		\right\|_{\ell^q}
		\asymp
		\|f\|_{M^{s,\theta}_{p,q}}.
	\]

	Finally,
	\[
		r_k^{d(1/p-1/2)}
		=
		|\xi_k|^{-\theta d(1/2-1/p)}
		\asymp
		\langle\xi_k\rangle^{-\theta d(1/2-1/p)}.
	\]
	Multiplying \eqref{eq:theta-sampling-block} by
	$w_\beta(\xi_k)=(1+|\xi_k|)^\beta$ and using
	$1+|\xi_k|\asymp\langle\xi_k\rangle$ gives
	\[
		w_\beta(\xi_k)
		\big\|(c_{j,k}(f))_j\big\|_{\ell^p}
		\asymp
		\langle\xi_k\rangle^s\|h_k\|_{L^p},
	\]
	because
	\[
		s=\beta-\theta d\left(\frac12-\frac1p\right).
	\]
	Taking the outer $\ell^q$ norm proves the asserted equivalence.
\end{proof}

\begin{remark}[The regime $\alpha>1$]
	If $\alpha>1$, then $\theta=1-\alpha<0$ and the frequency
	radii $r_k\asymp|\xi_k|^\theta$ tend to zero at high frequency.
	This is not a standard $\theta$--modulation covering, since the
	classical parameter range is $0\le\theta\le1$. Nevertheless,
	the family
	\[
	B_k=r_k Q+\xi_k
	\]
	is a structured admissible frequency covering: intersecting
	blocks have comparable $r_k$, and the overlap is uniformly
	bounded. Thus the corresponding coefficient spaces belong
	naturally to the general decomposition-space framework
	\cite{FeichtingerVoigtlaender2017,Voigtlaender2022}.
	The point specific to the present paper is the explicit
	phase-space realization of this shrinking-frequency geometry
	and its connection with the Beltr\'an--Bennett subdyadic scale.
\end{remark}

\subsection{Independence of the window and the lattice}

A priori, the above definition depends on the choice of window $\phi$ and lattice $\Omega_\alpha$. We will show that, within a reasonable class of windows and lattices adapted to $g^\alpha$, the norm $M^{p,q}_{\alpha,\beta}$ is well defined up to equivalence.

\begin{definition}[$\alpha$-admissible frame pairs]\label{def:alpha-adm}
	A pair $(\phi,\Omega_\alpha)$ is called $\alpha$-admissible if
	$\phi\in\mathcal S(\mathbb R^d)$ satisfies the window hypotheses of
	Proposition~\ref{prop:frame}, and $\Omega_\alpha$ is a product lattice
	obtained from the construction of Lemma~\ref{lem:covering} at an
	admissible scale $\rho_\ast$ for which Proposition~\ref{prop:frame}
	applies. In particular,
	\[
		\mathcal G_\alpha(\phi)
		=
		\{\phi_z:z\in\Omega_\alpha\}
	\]
	is a frame for $L^2(\mathbb R^d)$, its frequency blocks have uniformly
	finite overlap, and the associated canonical dual has the localization
	properties established in Proposition~\ref{prop:frame}.

	By a slight abuse of terminology, we call a window or a lattice
	$\alpha$-admissible when it occurs as part of an $\alpha$-admissible
	frame pair.
\end{definition}

\subsection{Existence of admissible pairs and norm equivalence}

\begin{proposition}[Existence of $\alpha$–admissible pairs]\label{prop:exist-adm}
	For every $\alpha>0$ there exist a family of subdyadic balls $\{B_k\}$ with finite overlap, centers $\{\xi_k\}$, and a nonzero window $\phi\in\mathcal{S}(\mathbb{R}^d)$ with compactly supported Fourier transform, bounded away from zero on a neighborhood of the origin, such that:
	\begin{enumerate}
		\item The packets $\phi_{x,\xi}$ are uniformly localized in $Q_\alpha(x,\xi;1)$.
		\item The lattice $\Omega_\alpha=\{(x_j^k,\xi_k)\}$ constructed in
Lemma~\ref{lem:covering} is relatively separated, and the
$d_{g^\alpha}$-balls centered in the region $|\xi|\ge1$ satisfy the
small-scale doubling estimate of Lemma~\ref{lem:volume}.
		\item $\mathcal{G}_\alpha(\phi)=\{\phi_{x_j^k,\xi_k}\}$ is a frame for $L^2(\mathbb{R}^d)$.
	\end{enumerate}
\end{proposition}

\begin{proof}
	The subdyadic partition with finite overlap is obtained by frequency coronas
	$|\xi|\sim2^k$ of thickness $\asymp2^{k(1-\alpha)}$ (cf.\ \cite{BeltranBennett2017}),
	which yields a family of $\alpha$–subdyadic balls $\{B_k\}$ in
	frequency with finite overlap and centers $\xi_k\in B_k$.
	
	If $\widehat\phi$ is supported in a fixed ball, then, by the definition of the dispersive packets
	\[
	\phi_{x,\xi}(t)
	=
	|\xi|^{\frac{d}{2}(1-\alpha)}\,\phi\bigl((t-x)|\xi|^{1-\alpha}\bigr)e^{it\cdot\xi},
	\]
	the Fourier transform $\widehat{\phi_{x,\xi}}$ is contained in a ball
	of radius $\asymp|\xi|^{1-\alpha}$ around $\xi$ by Lemma \ref{lem:freq-local-packets}, and $\phi_{x,\xi}$ is essentially
	concentrated in a neighborhood of $x$ of radius
	$\asymp|\xi|^{\alpha-1}$. This shows that $\phi_{x,\xi}$ is uniformly
	localized in $Q_\alpha(x,\xi;1)$ (in the sense of rapid decay outside
	that block), proving item (1).
	
	The geometric properties of the lattice $\Omega_\alpha=\{(x_j^k,\xi_k)\}$
	(covering, separation and relative separation) are contained in
	Lemmas~\ref{lem:covering} and \ref{lem:rel-sep-Omega}, while the frame
	property of the family $\mathcal{G}_\alpha(\phi)=\{\phi_{x_j^k,\xi_k}\}$
	was established in Proposition~\ref{prop:frame}. This gives (2) and (3) and
	completes the proof.
\end{proof}

\begin{lemma}[Moderateness of the radial weight]\label{lem:moderado}
	Let $w_\beta(\xi)=(1+|\xi|)^\beta$. There exists $C_\beta\ge1$ such that if
	\[
	\mathbf d_\alpha\big((x,\xi),(y,\eta)\big)\le1,
	\]
	then
	\[
	C_\beta^{-1}\,w_\beta(\xi)\ \le\ w_\beta(\eta)\ \le\ C_\beta\,w_\beta(\xi).
	\]
\end{lemma}

\begin{proof}
	Let
	\[
	r:=1+\min\{|\xi|,|\eta|\}.
	\]
	If
	\[
	\mathbf d_\alpha\big((x,\xi),(y,\eta)\big)\le1,
	\]
	then, by definition,
	\[
	r^{\alpha-1}|\xi-\eta|\le1,
	\]
	and hence
	\[
	|\xi-\eta|
	\le r^{1-\alpha}.
	\]
	Since $\alpha>0$ and $r\ge1$, we have
	\[
	r^{1-\alpha}\le r.
	\]
	Therefore
	\[
	\big||\xi|-|\eta|\big|
	\le|\xi-\eta|
	\le1+\min\{|\xi|,|\eta|\}.
	\]
	If
	\[
	m:=\min\{|\xi|,|\eta|\},
	\qquad
	M:=\max\{|\xi|,|\eta|\},
	\]
	then
	\[
	M\le m+(1+m)=1+2m,
	\]
	and consequently
	\[
	1+M\le2(1+m).
	\]
	Thus
	\[
	\frac12
	\le
	\frac{1+|\eta|}{1+|\xi|}
	\le2.
	\]
	Raising this estimate to the power $\beta$ gives
	\[
	2^{-|\beta|}(1+|\xi|)^\beta
	\le
	(1+|\eta|)^\beta
	\le
	2^{|\beta|}(1+|\xi|)^\beta.
	\]
	Hence the assertion holds with $C_\beta=2^{|\beta|}$.
\end{proof}

\begin{lemma}[Cross almost–diagonalization]\label{lem:cross-gram}
	Let $\phi,\psi$ be $\alpha$–admissible windows and $\Omega_\alpha,\Omega'_\alpha$ be $\alpha$–admissible lattices. Define
	\[
	G^{\psi,\phi}_{w,z}:=\langle \phi_w,\psi_z\rangle,\qquad  w\in\Omega_\alpha,\ z\in\Omega'_\alpha.
	\]
	Then, for every $N>0$, there exists $C_N<\infty$ such that
	\begin{equation}\label{eq:cross-gram-decay}
		|G^{\psi,\phi}_{w,z}|\ \le\ C_N\,\big(1+ \mathbf d_\alpha(w,z)\big)^{-N},
		\qquad w\in\Omega_\alpha,\ z\in\Omega'_\alpha.
	\end{equation}
	The constant $C_N$ depends only on $(d,\alpha)$, on the admissibility data of
$\Omega_\alpha,\Omega'_\alpha$, on the fixed Fourier-support radii of
$\phi,\psi$, and on a finite number of Schwartz seminorms of $\phi,\psi$, but
is independent of $w,z$.
\end{lemma}

\begin{proof}
	Write $w=(x,\xi)\in\Omega_\alpha$ and
	$z=(x',\xi')\in\Omega'_\alpha$, and set
	\[
		a:=|\xi|^{1-\alpha},
		\qquad
		b:=|\xi'|^{1-\alpha},
		\qquad
		r:=1+\min\{|\xi|,|\xi'|\}.
	\]
	Choose $R_\phi,R_\psi>0$ such that
	$\operatorname{supp}\widehat\phi\subset B(0,R_\phi)$ and
	$\operatorname{supp}\widehat\psi\subset B(0,R_\psi)$.
	By Lemma~\ref{lem:freq-local-packets}, if
	$\operatorname{supp}\widehat{\phi_w}$ and
	$\operatorname{supp}\widehat{\psi_z}$ are disjoint, then
	$G^{\psi,\phi}_{w,z}=0$, so there is nothing to prove.

	Assume therefore that the two Fourier supports intersect. Then
	\[
		|\xi-\xi'|
		\le R_\phi a+R_\psi b.
	\]
	Since $\alpha>0$, this overlap condition implies
	$|\xi|\asymp|\xi'|$ uniformly: for large frequencies this follows from
	$|\zeta|^{1-\alpha}=o(|\zeta|)$ as $|\zeta|\to\infty$, while the remaining
	bounded frequency region is absorbed into the constants. Consequently,
	\begin{equation}\label{eq:cross-comparable-scales}
		a\asymp b\asymp r^{1-\alpha},
		\qquad
		r^{\alpha-1}|\xi-\xi'|\lesssim1.
	\end{equation}

	Using Plancherel and the Fourier formula of
	Lemma~\ref{lem:freq-local-packets}, followed by the change of variables
	$\eta=\xi+a\zeta$, we obtain, up to the fixed Fourier-normalization
	constant,
	\[
	G^{\psi,\phi}_{w,z}
	=
	\Big(\frac{a}{b}\Big)^{d/2}
	e^{i x'\cdot(\xi-\xi')}
	\int_{\mathbb R^d}
	 e^{-i a(x-x')\cdot\zeta}\,
	 q_{w,z}(\zeta)\,d\zeta,
	\]
	where
	\[
	q_{w,z}(\zeta)
	:=
	\widehat\phi(\zeta)\,
	\overline{\widehat\psi\!\left(
		\frac{\xi-\xi'}{b}+\frac{a}{b}\zeta
	\right)}.
	\]
	By \eqref{eq:cross-comparable-scales}, the ratios $a/b$ and $b/a$ are
	uniformly bounded. Hence, for every integer $M\ge0$,
	\[
	\big\|(1-\Delta_\zeta)^M q_{w,z}\big\|_{L^1_\zeta}
	\le C_M,
	\]
	where $C_M$ is independent of $w,z$ and depends only on finitely many
	Schwartz seminorms of $\phi$ and $\psi$ (and on the fixed support radii).
	Integrating by parts in $\zeta$ therefore gives
	\[
	|G^{\psi,\phi}_{w,z}|
	\le
	C_M\bigl(1+a^2|x-x'|^2\bigr)^{-M}
	\le
	C_M\bigl(1+a|x-x'|\bigr)^{-2M}.
	\]
	Using again \eqref{eq:cross-comparable-scales}, for every $N>0$ we may
	choose $M$ large enough to obtain
	\[
	|G^{\psi,\phi}_{w,z}|
	\le
	C_N\bigl(1+r^{1-\alpha}|x-x'|\bigr)^{-N}.
	\]
	On the nonzero part of the cross Gramian the frequency component of the
adapted separation function is uniformly bounded by
	\eqref{eq:cross-comparable-scales}. Since
	\[
	\mathbf d_\alpha(w,z)
	\asymp
	r^{1-\alpha}|x-x'|
	+r^{\alpha-1}|\xi-\xi'|,
	\]
	the preceding estimate yields
	\[
	|G^{\psi,\phi}_{w,z}|
	\le C_N\bigl(1+\mathbf d_\alpha(w,z)\bigr)^{-N}.
	\]
	Together with the disjoint-support case, this proves
	\eqref{eq:cross-gram-decay}.
\end{proof}

\begin{theorem}[Window and lattice independence]\label{thm:window-grid-indep}
	Let $(\phi,\Omega_\alpha)$ and $(\psi,\Omega'_\alpha)$ be two
	subdyadic frame pairs obtained by the construction of
	Proposition~\ref{prop:frame}, with product lattices of the form
	constructed in Lemma~\ref{lem:covering}. For every
	$\beta\in\mathbb{R}$ there exist constants
	\[
	0<c_\beta\le C_\beta<\infty
	\]
	such that, simultaneously for all $1\le p,q\le\infty$,
	\[
	c_\beta\,
	\|f\|_{M^{p,q}_{\alpha,\beta}}^{(\phi,\Omega_\alpha)}
	\le
	\|f\|_{M^{p,q}_{\alpha,\beta}}^{(\psi,\Omega'_\alpha)}
	\le
	C_\beta\,
	\|f\|_{M^{p,q}_{\alpha,\beta}}^{(\phi,\Omega_\alpha)} .
	\]
	The constants may depend on $d,\alpha,\beta$, on the fixed window data,
	and on the geometric parameters of the two constructions, but they are
	independent of $p,q$ and $f$.
\end{theorem}

\begin{proof}
	Let $\{\widetilde\phi_z\}_{z\in\Omega_\alpha}$ be the canonical dual of
	the frame generated by $\phi$, and define
	\[
	H_{w,z}:=\langle\widetilde\phi_z,\psi_w\rangle,
	\qquad
	w\in\Omega'_\alpha,\quad z\in\Omega_\alpha.
	\]
	We first estimate $H$ directly. Write
	\[
	z=(x_j^k,\xi_k),
	\qquad
	w=(y_i^\ell,\eta_\ell),
	\]
	and set
	\[
	r_k:=|\xi_k|^{1-\alpha},
	\qquad
	s_\ell:=|\eta_\ell|^{1-\alpha}.
	\]
	By \eqref{eq:dual-fourier-packet}, the dual packet
	$\widetilde\phi_z$ has Fourier support contained in the same $k$-th
	frequency support as $\phi_z$, and its rescaled Fourier amplitude satisfies
	the uniform derivative bounds \eqref{eq:dual-amplitude-uniform}. Hence
	$H_{w,z}=0$ unless the $k$-th and $\ell$-th frequency supports intersect.
	For an interacting pair one has, uniformly,
	\begin{equation}\label{eq:window-grid-overlap-scales}
		r_k\asymp s_\ell,
		\qquad
		1+|\xi_k|\asymp1+|\eta_\ell|,
	\end{equation}
	and the frequency component of $\mathbf d_\alpha(w,z)$ is bounded.

	Repeating the Fourier-side rescaling and integration-by-parts argument of
	Lemma~\ref{lem:cross-gram}, with the uniformly controlled amplitude of
	$\widetilde\phi_z$ in place of $\widehat\phi$, gives, for every $N>0$,
	\begin{equation}\label{eq:dual-cross-direct}
		|H_{w,z}|
		\le
		C_N\bigl(1+\mathbf d_\alpha(w,z)\bigr)^{-N}.
	\end{equation}
	The constant is uniform in $w,z$.

	We next prove directly that $H$ is bounded on the mixed coefficient
	spaces. Fix an interacting frequency pair $(k,\ell)$. Since the spatial
	meshes of the two product lattices are respectively comparable to
	$r_k^{-1}$ and $s_\ell^{-1}$,
	\eqref{eq:window-grid-overlap-scales} and
	\eqref{eq:dual-cross-direct} imply, for $N>d$,
	\[
	\sup_i\sum_j |H_{(i,\ell),(j,k)}|
	+
	\sup_j\sum_i |H_{(i,\ell),(j,k)}|
	\le C,
	\]
	with $C$ independent of $k,\ell$. Thus every nonzero frequency block is
	bounded on $\ell^p(\mathbb Z^d)$, uniformly for
	$1\le p\le\infty$.

	By the uniformly finite overlap of the two frequency coverings together
	with \eqref{eq:window-grid-overlap-scales}, every frequency block
	interacts with only uniformly finitely many blocks of the other covering.
	On such interacting pairs,
	\[
	w_\beta(\eta_\ell)\asymp w_\beta(\xi_k).
	\]
	A finite-degree Schur estimate in the outer frequency index therefore
	yields
	\begin{equation}\label{eq:dual-cross-mixed-bound}
		\|Hc\|_{\ell^q(\ell^p;w_\beta)}
		\le
		C_\beta\|c\|_{\ell^q(\ell^p;w_\beta)},
		\qquad 1\le p,q\le\infty,
	\end{equation}
	with a constant independent of $p$ and $q$.

	For $f\in\mathcal S(\mathbb R^d)$, the frame reconstruction gives
	\[
	f
	=
	\sum_{z\in\Omega_\alpha}
	\langle f,\phi_z\rangle\widetilde\phi_z
	\quad\text{in }L^2,
	\]
	and hence
	\[
	C_\psi f=HC_\phi f.
	\]
	By \eqref{eq:dual-cross-mixed-bound},
	\[
	\|C_\psi f\|_{\ell^q(\ell^p;w_\beta)}
	\le
	C_\beta
	\|C_\phi f\|_{\ell^q(\ell^p;w_\beta)}.
	\]
	Interchanging $(\phi,\Omega_\alpha)$ and
	$(\psi,\Omega'_\alpha)$ gives the reverse inequality.

	For a general $f\in\mathcal S'(\mathbb R^d)$ with
$C_\phi f\in\ell^q(\ell^p;w_\beta)$, the localized frame expansion of each
test packet $\psi_w$,
\[
\psi_w
=
\sum_{z\in\Omega_\alpha}
\langle\psi_w,\widetilde\phi_z\rangle\,\phi_z,
\]
converges in $\mathcal S(\mathbb R^d)$. Pairing this identity with $f$
therefore gives
\[
C_\psi f=HC_\phi f.
\]
Hence \eqref{eq:dual-cross-mixed-bound} yields the same estimate for all
$1\le p,q\le\infty$. Interchanging the two frame pairs gives the reverse
inequality. This proves the claimed norm equivalence, with constants
uniform in $p$ and $q$.
\end{proof}

\noindent In particular, from now on we simply write
$M^{p,q}_{\alpha,\beta}$ without explicit reference to the window and the
lattice, with the understanding that the frame pair is chosen within the
class covered by Theorem~\ref{thm:window-grid-indep}.

\subsection{Relation with Sobolev spaces and density of
\texorpdfstring{$\mathcal{S}$}{S}}

\begin{proposition}\label{prop:L2}
	For every $\alpha>0$ and every $\beta\in\mathbb R$,
	\[
	M^{2,2}_{\alpha,\beta}
	=
	H^\beta(\mathbb R^d)
	\]
	with equivalence of norms. In particular,
	$M^{2,2}_{\alpha,0}=L^2(\mathbb R^d)$.
\end{proposition}

\begin{proof}
	Set
	\[
		r_k:=|\xi_k|^{1-\alpha},
		\qquad
		a_k:=b\rho_\ast r_k^{-1},
		\qquad
		\Phi_k(\eta):=
		\widehat\phi\!\left(\frac{\eta-\xi_k}{r_k}\right).
	\]
	We first justify the block Parseval identity at the level of tempered
	distributions. For fixed $k$, the compactly supported distribution
	\[
		F_k:=\widehat f\,\overline{\Phi_k}
	\]
	is supported in $B(\xi_k,R_1r_k)$. Since
	$a_kr_kR_1=b\rho_\ast R_1<\pi$, this support lies in the interior of a
	fundamental cube for the dual lattice $(2\pi/a_k)\mathbb Z^d$.
	By Lemma~\ref{lem:freq-local-packets}, the sequence
	$j\mapsto c_{j,k}(f)$ is, up to a fixed factor and a unimodular phase,
	the Fourier-coefficient sequence of the periodization of $F_k$ on that
	cube. Hence the standard Fourier-series characterization of $L^2$ on the
	torus implies
	\[
		(c_{j,k}(f))_{j\in\mathbb Z^d}\in\ell^2
		\quad\Longleftrightarrow\quad
		F_k\in L^2,
	\]
	and in this case Parseval gives
	\begin{equation}\label{eq:L2-block-parseval}
		\sum_{j\in\mathbb Z^d}
		|c_{j,k}(f)|^2
		=
		C_{\rho_\ast}
		\int_{\mathbb R^d}
		|\widehat f(\eta)|^2|\Phi_k(\eta)|^2\,d\eta,
	\end{equation}
	where the right-hand side is interpreted after identifying
	$F_k=\widehat f\,\overline{\Phi_k}$ with its $L^2$ representative, and
	$C_{\rho_\ast}>0$ is independent of $k$.

	Now suppose first that $f\in M^{2,2}_{\alpha,\beta}$. Then for every $k$
	the sequence $(c_{j,k}(f))_j$ belongs to $\ell^2$, so the preceding
	argument shows that $\widehat f\,\overline{\Phi_k}$ is an $L^2$ function.
	The lower bound for the function
	\[
		P(\eta):=\sum_k|\Phi_k(\eta)|^2
	\]
	proved in Proposition~\ref{prop:frame} implies that the open sets on
	which the $\Phi_k$ do not vanish cover frequency space. On each such set
	$\widehat f$ is therefore represented by an $L^2_{\mathrm{loc}}$ function;
	the local representatives agree on overlaps because they represent the
	same tempered distribution. Thus $\widehat f$ itself may be identified
	with an $L^2_{\mathrm{loc}}$ function.

	Multiplying \eqref{eq:L2-block-parseval} by
	$w_\beta(\xi_k)^2=(1+|\xi_k|)^{2\beta}$ and summing over $k$ gives
	\begin{equation}\label{eq:L2-weighted-parseval}
		\|f\|_{M^{2,2}_{\alpha,\beta}}^2
		=
		C_{\rho_\ast}
		\int_{\mathbb R^d}
		|\widehat f(\eta)|^2P_\beta(\eta)\,d\eta,
	\end{equation}
	where
	\[
		P_\beta(\eta)
		:=
		\sum_k
		(1+|\xi_k|)^{2\beta}|\Phi_k(\eta)|^2.
	\]
	Whenever $\Phi_k(\eta)\neq0$, the frequency localization gives
	\[
		1+|\xi_k|\asymp1+|\eta|,
	\]
	with constants independent of $k$ and $\eta$. Consequently,
	\[
		P_\beta(\eta)
		\asymp
		(1+|\eta|)^{2\beta}P(\eta)
		\asymp
		\langle\eta\rangle^{2\beta},
	\]
	because $P$ is bounded above and below by positive constants and
	$1+|\eta|\asymp\langle\eta\rangle$. Equation
	\eqref{eq:L2-weighted-parseval} therefore yields
	\[
		\|f\|_{M^{2,2}_{\alpha,\beta}}
		\asymp
		\|\langle\eta\rangle^\beta\widehat f(\eta)\|_{L^2_\eta}
		\asymp
		\|f\|_{H^\beta}.
	\]
	In particular, $f\in H^\beta$.

	Conversely, let $f\in H^\beta(\mathbb R^d)$. Then each
	$F_k=\widehat f\,\overline{\Phi_k}$ belongs to $L^2$, so the same
	Fourier-series Parseval identity \eqref{eq:L2-block-parseval} applies.
	Summing it with the weights as above and using
	$P_\beta(\eta)\asymp\langle\eta\rangle^{2\beta}$ gives
	\[
		\|f\|_{M^{2,2}_{\alpha,\beta}}
		\lesssim
		\|f\|_{H^\beta}.
	\]
	Thus $H^\beta\subset M^{2,2}_{\alpha,\beta}$, and the reverse inclusion
	and equivalence of norms were proved above. The case $\beta=0$ gives
	$M^{2,2}_{\alpha,0}=L^2$.
\end{proof}

\begin{proposition}\label{prop:densidad}
	If $1\le p,q<\infty$, then $\mathcal{S}(\mathbb{R}^d)$ is dense in $M^{p,q}_{\alpha,\beta}$.
\end{proposition}

\begin{proof}
	Set
	\[
	X:=\ell^q(\ell^p;w_\beta),
	\qquad
	C_\phi f:=(\langle f,\phi_z\rangle)_{z\in\Omega_\alpha}.
	\]
	We first record the synthesis estimate that will also justify the
	Banach-space realization of the coefficient norm. For $c\in X$, define
	\[
	S_\phi c:=\sum_{z\in\Omega_\alpha}c_z\phi_z
	\]
	as a weak sum in $\mathcal S'(\mathbb R^d)$. This sum is well defined:
	for every Schwartz test function, its coefficients against the packets
	$\phi_z$ decrease faster than any polynomial in the phase-space index,
	whereas an element of $X$ has at most the growth allowed by the fixed
	radial weight $w_\beta$.

	The coefficients of $S_\phi c$ are
\[
(C_\phi S_\phi c)_w
=
\sum_z c_z\langle\phi_z,\phi_w\rangle.
\]
We use here the same block--Schur mechanism as in the proof of
Theorem~\ref{thm:window-grid-indep}. The primal Gram matrix vanishes
between disjoint frequency blocks. On interacting blocks the frequency
scales are uniformly comparable, the weights $w_\sigma$ are comparable,
and the off--diagonal estimate of Proposition~\ref{prop:frame} gives,
for every $N>d$, uniform row and column sums in the spatial lattice
indices. Since the frequency-overlap graph has uniformly finite degree,
we obtain, for every $\sigma\in\mathbb R$ and every
$1\le p,q\le\infty$,
\begin{equation}\label{eq:primal-gram-block-bounded}
G_\phi:
\ell^q(\ell^p;w_\sigma)
\longrightarrow
\ell^q(\ell^p;w_\sigma)
\quad\text{boundedly},
\end{equation}
where $(G_\phi)_{w,z}=\langle\phi_z,\phi_w\rangle$. In particular,
with $\sigma=\beta$,
\[
\|S_\phi c\|_{M^{p,q}_{\alpha,\beta}}
=
\|C_\phi S_\phi c\|_X
\le C\|c\|_X.
\]
Thus
\[
S_\phi:X\longrightarrow M^{p,q}_{\alpha,\beta}
\]
is bounded.

Let $\{\widetilde\phi_z\}_{z\in\Omega_\alpha}$ be the canonical dual
frame and set
\[
D_{w,z}:=\langle\widetilde\phi_z,\widetilde\phi_w\rangle.
\]
By \eqref{eq:dual-fourier-packet}, the dual packets have the same
frequency supports as the corresponding primal packets, and
\eqref{eq:dual-gram-direct} gives arbitrary off--diagonal decay on the
interacting blocks. The identical block--Schur argument therefore gives,
for every $\sigma\in\mathbb R$ and every $1\le p,q\le\infty$,
\begin{equation}\label{eq:dual-gram-block-bounded}
D:
\ell^q(\ell^p;w_\sigma)
\longrightarrow
\ell^q(\ell^p;w_\sigma)
\quad\text{boundedly}.
\end{equation}
	For $f\in M^{p,q}_{\alpha,\beta}$, write $c=C_\phi f\in X$.
For every $h\in\mathcal S(\mathbb R^d)$, the localized dual-frame
expansion
\[
h
=
\sum_z
\langle h,\widetilde\phi_z\rangle\,\phi_z
\]
converges in $\mathcal S(\mathbb R^d)$. Indeed, the coefficients of a
Schwartz function against the packets have rapid phase-space decay,
the dual Gram matrix $D$ has arbitrarily high off--diagonal decay by
\eqref{eq:dual-gram-direct}, and the Schwartz seminorms of the packets
grow at most polynomially in the phase-space indices.

Since
\[
\big(\langle h,\widetilde\phi_z\rangle\big)_z
=
D\,C_\phi h
\]
and $D$ is self-adjoint, being the Gram matrix of the canonical dual
frame. Pairing the preceding expansion with $f$ gives
\[
\begin{aligned}
\langle f,h\rangle
&=
\sum_z c_z\,
\overline{\langle h,\widetilde\phi_z\rangle}  \\
&=
\sum_w (Dc)_w\,
\overline{\langle h,\phi_w\rangle} \\
&=
\langle S_\phi Dc,h\rangle .
\end{aligned}
\]
Hence
\begin{equation}\label{eq:analysis-synthesis-retraction}
f=S_\phi D C_\phi f
\qquad\text{in }\mathcal S'(\mathbb R^d).
\end{equation}
	The rearrangement is justified by the block bounds above together with
the rapid decay of the test-function coefficients. In particular,
$S_\phi D$ is a bounded	left inverse of $C_\phi$ on $M^{p,q}_{\alpha,\beta}$. Hence
	$C_\phi(M^{p,q}_{\alpha,\beta})$ is the range of the bounded projection
	$C_\phi S_\phi D$ on $X$ and is therefore closed. Since the norm of
	$M^{p,q}_{\alpha,\beta}$ is the pullback of the norm of $X$ under
	$C_\phi$, this also proves completeness. The same argument, without the
	finite-sequence approximation below, applies at the endpoint values
	$p=\infty$ or $q=\infty$.

	Now assume $1\le p,q<\infty$. Set $d:=DC_\phi f\in X$. Since finitely
	supported sequences are dense in $X$, if $d^F$ denotes the restriction of
	$d$ to a finite set $F\subset\Omega_\alpha$, then $d^F\to d$ in $X$. Define
	\[
	f_F:=S_\phi d^F
	=\sum_{z\in F}d_z\phi_z.
	\]
	Every $f_F$ belongs to $\mathcal S(\mathbb R^d)$, and by
	\eqref{eq:analysis-synthesis-retraction} and the synthesis estimate,
	\[
	\|f-f_F\|_{M^{p,q}_{\alpha,\beta}}
	=
	\|S_\phi(d-d^F)\|_{M^{p,q}_{\alpha,\beta}}
	\le
	C\|d-d^F\|_X
	\longrightarrow0.
	\]
	Thus $\mathcal S(\mathbb R^d)$ is dense in
	$M^{p,q}_{\alpha,\beta}$.
\end{proof}

\subsection{Duality and basic inclusions}

\begin{proposition}[Duality]\label{prop:dualidad}
	Let $1<p,q<\infty$ and $\beta\in\mathbb{R}$. Then
	\[
	(M^{p,q}_{\alpha,\beta})' \simeq M^{p',q'}_{\alpha,-\beta},
	\]
	with the pairing obtained by continuous extension of the $L^2$ pairing
	\[
	\langle f,g\rangle
	=
	\int_{\mathbb R^d}f(x)\,\overline{g(x)}\,dx,
	\qquad f,g\in\mathcal S(\mathbb R^d).
	\]
\end{proposition}

\begin{proof}
	Let $\{\widetilde\phi_z\}_{z\in\Omega_\alpha}$ be the canonical dual
	frame and write
	\[
	c_z(f):=\langle f,\phi_z\rangle,
	\qquad
	\widetilde c_z(f):=\langle f,\widetilde\phi_z\rangle
	\quad (f\in\mathcal S).
	\]
	For $f,g\in\mathcal S(\mathbb R^d)$, the $L^2$ frame identity gives
	\[
	\langle f,g\rangle
	=
	\sum_{z\in\Omega_\alpha}
	c_z(f)\,\overline{\widetilde c_z(g)}.
	\]
	If
	\[
	D_{w,z}:=\langle\widetilde\phi_z,\widetilde\phi_w\rangle,
	\]
	then the dual-frame expansion implies
$\widetilde c(g)=D\,c(g)$ for $g\in\mathcal S$. By
\eqref{eq:dual-gram-block-bounded}, applied with $\sigma=-\beta$ and
exponents $p',q'$, we obtain
\[
\|\widetilde c(g)\|_{\ell^{q'}(\ell^{p'};w_{-\beta})}
\lesssim
\|g\|_{M^{p',q'}_{\alpha,-\beta}}.
\]
	Discrete H\"older's inequality now yields
	\[
	|\langle f,g\rangle|
	\lesssim
	\|f\|_{M^{p,q}_{\alpha,\beta}}
	\|g\|_{M^{p',q'}_{\alpha,-\beta}},
	\qquad f,g\in\mathcal S.
	\]
	By Proposition~\ref{prop:densidad}, this pairing extends continuously to
	$M^{p,q}_{\alpha,\beta}\times M^{p',q'}_{\alpha,-\beta}$. Hence every
	$g\in M^{p',q'}_{\alpha,-\beta}$ defines a bounded functional on
	$M^{p,q}_{\alpha,\beta}$.

	Conversely, let $\Lambda\in(M^{p,q}_{\alpha,\beta})'$ and set
	\[
	X:=\ell^q(\ell^p;w_\beta),
	\qquad
	Y:=C_\phi(M^{p,q}_{\alpha,\beta})\subset X.
	\]
	Since the norm of $M^{p,q}_{\alpha,\beta}$ is the coefficient norm,
	\[
	\lambda(C_\phi f):=\Lambda(f),
	\qquad C_\phi f\in Y,
	\]
	defines a bounded functional on the subspace $Y$ of $X$. By the
	Hahn--Banach theorem, $\lambda$ extends to a bounded functional on $X$.
	Because $1<p,q<\infty$,
	\[
	X'
	=
	\ell^{q'}(\ell^{p'};w_{-\beta}),
	\]
	so there exists $d\in\ell^{q'}(\ell^{p'};w_{-\beta})$ such that
	\[
	\Lambda(f)
	=
	\sum_{z\in\Omega_\alpha}c_z(f)\,\overline{d_z},
	\qquad f\in M^{p,q}_{\alpha,\beta},
	\]
	with $\|d\|\lesssim\|\Lambda\|$.

	Define, using the bounded primal synthesis operator constructed in the
	proof of Proposition~\ref{prop:densidad},
	\[
	g:=S_\phi d
	=
	\sum_{z\in\Omega_\alpha}d_z\phi_z.
	\]
	The same synthesis estimate with the exponents $p',q'$ and weight
	$w_{-\beta}$ gives
	\[
	g\in M^{p',q'}_{\alpha,-\beta},
	\qquad
	\|g\|_{M^{p',q'}_{\alpha,-\beta}}
	\lesssim\|d\|.
	\]
	For $f\in\mathcal S(\mathbb R^d)$, finite truncation of the synthesis
	series gives
	\[
	\langle f,g\rangle
	=
	\sum_z c_z(f)\,\overline{d_z}
	=
	\Lambda(f).
	\]
	Both sides are continuous on $M^{p,q}_{\alpha,\beta}$, and
	Proposition~\ref{prop:densidad} therefore extends the identity to every
	$f\in M^{p,q}_{\alpha,\beta}$.

	If two elements of $M^{p',q'}_{\alpha,-\beta}$ represent the same
	functional, their difference pairs to zero with every Schwartz function
	and hence vanishes in $\mathcal S'$. Thus the representing element is
	unique. The preceding estimates in both directions give equivalence of
	norms and complete the proof.
\end{proof}

\begin{proposition}[Monotone inclusions]\label{prop:inclusiones}
	If $1\le p_1\le p_2\le\infty$ and $1\le q_1\le q_2\le\infty$, then
	\[
	M^{p_1,q_1}_{\alpha,\beta}\hookrightarrow M^{p_2,q_2}_{\alpha,\beta}
	\]
	continuously.
\end{proposition}

\begin{proof}
	The discrete inclusions $\ell^{p_1}\subset \ell^{p_2}$ for $p_1\le p_2$ (with
	$\|\cdot\|_{\ell^{p_2}}\le \|\cdot\|_{\ell^{p_1}}$) and similarly
	$\ell^{q_1}\subset \ell^{q_2}$ imply
	\[
	\|c(f)\|_{\ell^{q_2}(\ell^{p_2};w_\beta)}\ \le\ \|c(f)\|_{\ell^{q_1}(\ell^{p_1};w_\beta)}.
	\]
	By the definition of the $M^{p,q}_{\alpha,\beta}$ norms via Gabor coefficients,
	this yields
	\[
	\|f\|_{M^{p_2,q_2}_{\alpha,\beta}}
	\lesssim
	\|f\|_{M^{p_1,q_1}_{\alpha,\beta}},
	\]
	which proves the continuous embedding.
\end{proof}

\medskip

\noindent
These properties show that $M^{p,q}_{\alpha,\beta}$ provides a robust and flexible
framework to measure regularity and mass distribution in the subdyadic geometry
induced by $g^\alpha$.

\section{Miyachi multipliers on dispersive modulation spaces}\label{seccion:Multiplicadores:Miyachi}

In this section we study Fourier multipliers adapted to the subdyadic
frequency scale. We use the high-frequency H\"ormander--Miyachi condition
\eqref{Cond:H-Miya:0}, which is the condition naturally associated with the
subdyadic balls of Beltr\'an and Bennett \cite{BeltranBennett2017}.

To avoid confusion with the weight parameter $\beta$ in
$M^{p,q}_{\alpha,\beta}$, we denote by $\mu\in\mathbb R$ the order of the
multiplier. Thus, at high frequency, the pointwise Miyachi condition reads
\[
|D^\gamma m(\xi)|
\lesssim
|\xi|^{-\mu+(\alpha-1)|\gamma|},
\qquad |\xi|^\alpha\ge1,
\]
for $|\gamma|\le N_d$, where
\[
N_d=\left\lfloor\frac d2\right\rfloor+1.
\]
The weaker local quadratic condition \eqref{Cond:H-Miya:0}, with $\beta$
there replaced by $\mu$, will be sufficient for the arguments below.
We assume throughout this section that the multiplier is supported in the
high-frequency region $\{|\xi|^\alpha\ge1\}$. A separate smooth
low-frequency component can be treated independently.

The point of the argument is not a general multiplier principle for
decomposition spaces, for which an extensive theory already exists, but an
explicit coefficient estimate for packets adapted to the
Beltr\'an--Bennett subdyadic geometry. The resulting matrix estimate will
yield boundedness on the full mixed-norm scale
$M^{p,q}_{\alpha,\beta}$.

\subsection{Microlocal set-up}

Let $m:\mathbb{R}^d\setminus\{0\}\to\mathbb{C}$ be supported in
$\{|\xi|^\alpha\ge1\}$ and satisfy the H\"ormander--Miyachi condition
\eqref{Cond:H-Miya:0} with parameters $(\alpha,\mu)$ at high frequency.
Initially, for $f\in\mathcal S(\mathbb R^d)$, we consider the Fourier
multiplier
\[
\widehat{T_m f}(\xi)=m(\xi)\widehat f(\xi).
\]
A bounded realization on the full modulation-space scale will be constructed
in Theorem~\ref{thm:Miyachi-M22}. In time--frequency terms, $T_m$ is a
pseudodifferential operator with symbol depending only on $\xi$. It preserves
the Fourier support of each packet, although it may modify its spatial
profile.

Let $\phi$ be an $\alpha$--admissible window and $\Omega_\alpha$ a subdyadic
Gabor lattice as in the previous sections. For $f\in\mathcal S(\mathbb R^d)$,
the Gabor coefficients of $T_m f$ can be written as
\[
c_{j,k}(T_m f)=\langle T_m f,\phi_{x_j^k,\xi_k}\rangle
=(2\pi)^{-d}\int_{\mathbb{R}^d}
m(\xi)\,\widehat{f}(\xi)\,
\overline{\widehat{\phi_{x_j^k,\xi_k}}(\xi)}\,d\xi.
\]
Since $\widehat{\phi_{x_j^k,\xi_k}}$ is supported on an
$\alpha$--subdyadic block around $\xi_k$
(Section~\ref{seccion:geoSubD:fase-espacio}), only the local behavior of
$m$ on that block enters the matrix coefficients. The next lemma records
the scale-invariant consequence of the H\"ormander--Miyachi condition that
will be used in the block estimate.

\begin{lemma}[Rescaled Miyachi control]\label{lem:m-local}
	Let $B_k^\ast$ be a fixed enlargement of a high-frequency
	$\alpha$--subdyadic block $B_k$, with centre $\xi_k$ and radius
	\[
	r_k\asymp |\xi_k|^{1-\alpha}.
	\]
	If $m$ satisfies \eqref{Cond:H-Miya:0} with order $\mu$, then for every
	multi-index $|\gamma|\le N_d$,
	\[
	r_k^{|\gamma|-d/2}
	\|D^\gamma m\|_{L^2(B_k^\ast)}
	\lesssim
	|\xi_k|^{-\mu},
	\]
	uniformly in $k$.
\end{lemma}

\begin{proof}
	Since $B_k^\ast$ is again an $\alpha$--subdyadic ball up to fixed
	constants,
	\[
	|B_k^\ast|\asymp r_k^d,
	\qquad
	\operatorname{dist}(B_k^\ast,0)\asymp|\xi_k|.
	\]
	The H\"ormander--Miyachi condition gives
	\[
	|B_k^\ast|^{-1/2}
	\|D^\gamma m\|_{L^2(B_k^\ast)}
	\lesssim
	|\xi_k|^{-\mu+(\alpha-1)|\gamma|}.
	\]
	Therefore
	\[
	r_k^{|\gamma|-d/2}
	\|D^\gamma m\|_{L^2(B_k^\ast)}
	\lesssim
	r_k^{|\gamma|}
	|\xi_k|^{-\mu+(\alpha-1)|\gamma|}.
	\]
	Using $r_k\asymp|\xi_k|^{1-\alpha}$, the scale factors cancel:
	\[
	r_k^{|\gamma|}
	|\xi_k|^{(\alpha-1)|\gamma|}
	\asymp1.
	\]
	Hence
	\[
	r_k^{|\gamma|-d/2}
	\|D^\gamma m\|_{L^2(B_k^\ast)}
	\lesssim|\xi_k|^{-\mu},
	\]
	as claimed.
\end{proof}

\begin{lemma}[Dual localization and Moore--Penrose identities]\label{lem:dual-local}
	Let $\{\phi_z\}_{z\in\Omega_\alpha}$ be a subdyadic Gabor frame of the
	form constructed in Proposition~\ref{prop:frame}, and let
	$\{\widetilde\phi_z\}_{z\in\Omega_\alpha}$ be its canonical dual. Then,
	for every $N>0$,
	\begin{equation}\label{eq:dual-gram-direct}
		|\langle\widetilde\phi_z,\widetilde\phi_w\rangle|
		+
		|\langle\widetilde\phi_z,\phi_w\rangle|
		\le
		C_N\bigl(1+\mathbf d_\alpha(w,z)\bigr)^{-N},
	\end{equation}
	with $C_N$ independent of $w,z$.

	Moreover, if
	\[
	C:L^2(\mathbb R^d)\to\ell^2(\Omega_\alpha),
	\qquad
	Cf=(\langle f,\phi_w\rangle)_{w\in\Omega_\alpha},
	\]
	$S=C^\ast C$, and $G=CC^\ast$, then the Moore--Penrose pseudoinverse of
	$G$ satisfies
	\begin{equation}\label{eq:gram-pseudoinverse-identities}
	(G^\dagger)_{w,z}
	=
	\langle\widetilde\phi_z,\widetilde\phi_w\rangle,
	\qquad
	(GG^\dagger)_{w,z}
	=
	\langle\widetilde\phi_z,\phi_w\rangle.
	\end{equation}
	In particular,
	\[
	GG^\dagger=CS^{-1}C^\ast
	\]
	is the orthogonal projection of $\ell^2(\Omega_\alpha)$ onto
	$\operatorname{ran}C$.
\end{lemma}

\begin{proof}
	The localization estimate \eqref{eq:dual-gram-direct} is precisely the
	direct Fourier-side localization proved in
	Proposition~\ref{prop:frame}; see
	\eqref{eq:dual-fourier-packet}--\eqref{eq:dual-direct-localization}, with
	the packet notation $\psi_z$ there replaced by $\phi_z$ here. Thus no
	inverse-closedness theorem for a global matrix algebra is needed.

	It remains to record the operator identities. Since $\{\phi_z\}$ is a
	frame, $C$ has closed range and $S=C^\ast C$ is boundedly invertible on
	$L^2$. On $\operatorname{ran}C$ one has
	\[
	G(Cf)=CC^\ast Cf=CSf,
	\]
	whereas $G=0$ on $(\operatorname{ran}C)^\perp$. Hence
	\[
	G^\dagger=CS^{-2}C^\ast.
	\]
	With the convention
	$\langle f,g\rangle=\int f\,\overline g$, we have $C^\ast e_z=\phi_z$,
	and therefore
	\[
	(G^\dagger)_{w,z}
	=
	\langle S^{-2}\phi_z,\phi_w\rangle
	=
	\langle S^{-1}\phi_z,S^{-1}\phi_w\rangle
	=
	\langle\widetilde\phi_z,\widetilde\phi_w\rangle.
	\]
	Furthermore,
	\[
	GG^\dagger
	=
	CC^\ast CS^{-2}C^\ast
	=
	CS^{-1}C^\ast,
	\]
	which is the orthogonal projection onto $\operatorname{ran}C$, and
	\[
	(GG^\dagger)_{w,z}
	=
	\langle S^{-1}\phi_z,\phi_w\rangle
	=
	\langle\widetilde\phi_z,\phi_w\rangle.
	\]
	This proves the assertions.
\end{proof}

\begin{lemma}[Block almost diagonalization of $T_m$]
\label{lem:matrix-Tm}
	For
	\[
	\mathcal G^m_{(j',k'),(j,k)}
	:=
	\big\langle
	T_m\phi_{x_j^k,\xi_k},
	\phi_{x_{j'}^{k'},\xi_{k'}}
	\big\rangle,
	\]
	the following properties hold.

	\begin{enumerate}
		\item If the Fourier supports of the $k$-th and $k'$-th packet
		families do not intersect, then
		\[
		\mathcal G^m_{(j',k'),(j,k)}=0.
		\]
		Consequently each frequency index interacts with only uniformly
		finitely many other frequency indices.

		\item Whenever the two frequency supports intersect,
		\[
		\sup_{j\in\mathbb Z^d}
		\sum_{j'\in\mathbb Z^d}
		\big|
		\mathcal G^m_{(j',k'),(j,k)}
		\big|
		+
		\sup_{j'\in\mathbb Z^d}
		\sum_{j\in\mathbb Z^d}
		\big|
		\mathcal G^m_{(j',k'),(j,k)}
		\big|
		\lesssim
		(1+|\xi_k|)^{-\mu},
		\]
		with a constant independent of $k,k'$.

		\item Consequently, for every
		$1\le p,q\le\infty$ and every $\beta\in\mathbb R$,
		\[
		\mathcal G^m:
		\ell^q(\ell^p;w_{\beta-\mu})
		\longrightarrow
		\ell^q(\ell^p;w_\beta)
		\]
		is bounded.
	\end{enumerate}
\end{lemma}

\begin{proof}
	The Fourier transform of a packet has the form
	\[
	\widehat{\phi_{x_j^k,\xi_k}}(\xi)
	=
	r_k^{-d/2}
	e^{-ix_j^k\cdot(\xi-\xi_k)}
	\widehat\phi
	\left(
	\frac{\xi-\xi_k}{r_k}
	\right),
	\qquad
	r_k=|\xi_k|^{1-\alpha}.
	\]
	Hence the matrix coefficient vanishes unless the corresponding
	frequency supports intersect. By the finite-overlap property of the
	covering, for each $k$ there are only uniformly finitely many such
	$k'$. For intersecting blocks one has
	\[
	r_{k'}\asymp r_k,
	\qquad
	1+|\xi_{k'}|\asymp1+|\xi_k|.
	\]

	Fix such a pair $(k,k')$ and put
	\[
	\lambda_{k,k'}:=\frac{r_k}{r_{k'}},
	\qquad
	\delta_{k,k'}:=\frac{\xi_k-\xi_{k'}}{r_{k'}}.
	\]
	Both $\lambda_{k,k'}$ and $\delta_{k,k'}$ remain in fixed compact
	sets, with $\lambda_{k,k'}$ bounded away from zero.

	After the change of variables
	\[
	\xi=\xi_k+r_k\eta,
	\]
	the coefficient can be written, up to a unimodular constant, as
	\[
	\mathcal G^m_{(j',k'),(j,k)}
	=
	\int_{\mathbb R^d}
	a_{k,k'}(\eta)
	e^{-i r_k(x_j^k-x_{j'}^{k'})\cdot\eta}
	\,d\eta,
	\]
	where
	\[
	a_{k,k'}(\eta)
	=
	\lambda_{k,k'}^{d/2}
	m(\xi_k+r_k\eta)
	\widehat\phi(\eta)
	\overline{
	\widehat\phi(
	\delta_{k,k'}+\lambda_{k,k'}\eta
	)}.
	\]
	The functions $a_{k,k'}$ are supported in one fixed compact set,
	independent of $k,k'$. By Leibniz' rule and
	Lemma~\ref{lem:m-local},
	\[
	\sum_{|\gamma|\le N_d}
	\|\partial^\gamma a_{k,k'}\|_{L^2}
	\lesssim
	(1+|\xi_k|)^{-\mu}.
	\]

	We use the following elementary discrete Sobolev estimate. If
	$a$ is supported in a fixed compact set, $N_d>d/2$, and
	$Y\subset\mathbb R^d$ is uniformly separated, then
	\begin{equation}\label{eq:discrete-sobolev-sampling}
		\sum_{y\in Y}|\check a(y)|
		\lesssim
		\sum_{|\gamma|\le N_d}
		\|\partial^\gamma a\|_{L^2}.
	\end{equation}
	Indeed, since $\check a$ and
	$y^\gamma\check a(y)$ are bandlimited to a fixed compact set, the
	standard separated-sampling estimate gives
	\[
	\sum_{y\in Y}
	|y^\gamma\check a(y)|^2
	\lesssim
	\|\partial^\gamma a\|_{L^2}^2.
	\]
	Cauchy--Schwarz, together with
	\[
	\sum_{y\in Y}(1+|y|)^{-2N_d}<\infty
	\qquad(N_d>d/2),
	\]
	yields \eqref{eq:discrete-sobolev-sampling}.

	Now
	\[
	x_j^k=b\rho_\ast r_k^{-1}j,
	\]
	and therefore
	\[
	r_k(x_j^k-x_{j'}^{k'})
	=
	b\rho_\ast
	\left(
	j-\lambda_{k,k'}j'
	\right).
	\]
	For fixed $j$, the set
	\[
	\left\{
	b\rho_\ast
	(j-\lambda_{k,k'}j'):
	j'\in\mathbb Z^d
	\right\}
	\]
	is uniformly separated, and the same is true after fixing $j'$ and
	varying $j$. Applying
	\eqref{eq:discrete-sobolev-sampling} to $a_{k,k'}$ gives
	\[
	\sup_j\sum_{j'}
	|\mathcal G^m_{(j',k'),(j,k)}|
	+
	\sup_{j'}\sum_j
	|\mathcal G^m_{(j',k'),(j,k)}|
	\lesssim
	(1+|\xi_k|)^{-\mu}.
	\]
	This proves the block Schur bounds.

	By the discrete Schur test, each nonzero $(k',k)$ block therefore
	defines an operator on $\ell^p(\mathbb Z^d)$, uniformly for
	$1\le p\le\infty$, with norm
	\[
	\lesssim(1+|\xi_k|)^{-\mu}.
	\]
	For interacting frequency blocks,
	\[
	1+|\xi_{k'}|\asymp1+|\xi_k|,
	\]
	and hence
	\[
	w_\beta(\xi_{k'})
	(1+|\xi_k|)^{-\mu}
	\lesssim
	w_{\beta-\mu}(\xi_k).
	\]
	Finally, the frequency-interaction graph has uniformly finite degree.
	Applying the finite-degree Schur estimate in the outer $k$ variable
	gives
	\[
	\|\mathcal G^m c\|_{\ell^q(\ell^p;w_\beta)}
	\lesssim
	\|c\|_{\ell^q(\ell^p;w_{\beta-\mu})}
	\]
	for every $1\le p,q\le\infty$.
\end{proof}

\subsection{Boundedness on the mixed-norm scale}

\begin{theorem}[Miyachi multipliers on dispersive modulation spaces]
\label{thm:Miyachi-M22}
	Let $m$ satisfy the hypotheses of
	Lemma~\ref{lem:matrix-Tm} with order $\mu\in\mathbb R$. Then, for every
	\[
	1\le p,q\le\infty,
	\qquad
	\beta\in\mathbb R,
	\]
	the Fourier multiplier initially defined on $\mathcal S(\mathbb R^d)$
	admits a bounded extension
	\[
	T_m:
	M^{p,q}_{\alpha,\beta-\mu}
	\longrightarrow
	M^{p,q}_{\alpha,\beta},
	\]
	and
	\[
	\|T_m f\|_{M^{p,q}_{\alpha,\beta}}
	\lesssim
	\|f\|_{M^{p,q}_{\alpha,\beta-\mu}}.
	\]
	If $p,q<\infty$, this extension is unique. At the endpoint values it is
	understood as the bounded coefficient-space realization constructed in
	the proof below. In particular, if $\mu=0$, then
	\[
	T_m:
	M^{p,q}_{\alpha,\beta}
	\longrightarrow
	M^{p,q}_{\alpha,\beta}
	\]
	is bounded for every $1\le p,q\le\infty$.
\end{theorem}

\begin{proof}
	Fix the primal frame $\{\phi_w\}_{w\in\Omega_\alpha}$ and let
	$\{\widetilde\phi_w\}$ be its canonical dual. Put
	\[
	X_\sigma:=\ell^q(\ell^p;w_\sigma),
	\qquad
	C_\phi f=(\langle f,\phi_w\rangle)_w,
	\qquad
	D:=G^\dagger.
	\]
	By \eqref{eq:dual-gram-block-bounded},
	\[
	D:X_\sigma\longrightarrow X_\sigma
	\]
	boundedly for every $\sigma\in\mathbb R$ and all
	$1\le p,q\le\infty$. Moreover, the proof of
	Proposition~\ref{prop:densidad} gives a bounded primal synthesis operator
	$S_\phi:X_\sigma\to M^{p,q}_{\alpha,\sigma}$ and the retraction identity
	\[
	R_\phi C_\phi=I,
	\qquad
	R_\phi:=S_\phi D,
	\]
	on $M^{p,q}_{\alpha,\sigma}$; see
	\eqref{eq:analysis-synthesis-retraction}.

	For $f\in M^{p,q}_{\alpha,\beta-\mu}$ we define
	\begin{equation}\label{eq:Tm-coefficient-extension}
	\mathcal T_m f
	:=
	R_\phi\,\mathcal G^m D C_\phi f.
	\end{equation}
	Lemma~\ref{lem:matrix-Tm}, together with the boundedness of $D$ and
	$R_\phi$, gives
	\[
	\begin{aligned}
	\|\mathcal T_m f\|_{M^{p,q}_{\alpha,\beta}}
	&\lesssim
	\|\mathcal G^m D C_\phi f\|_{X_\beta} \\
	&\lesssim
	\|D C_\phi f\|_{X_{\beta-\mu}} \\
	&\lesssim
	\|C_\phi f\|_{X_{\beta-\mu}}
	=
	\|f\|_{M^{p,q}_{\alpha,\beta-\mu}}.
	\end{aligned}
	\]
	Thus \eqref{eq:Tm-coefficient-extension} defines a bounded operator on the
	full mixed-norm range, including the endpoints.

	It remains to verify that this operator agrees with the Fourier multiplier
	on test functions. Let $f\in\mathcal S(\mathbb R^d)$ and set
	\[
	a:=D C_\phi f
	=
	(\langle f,\widetilde\phi_r\rangle)_r.
	\]
	The localized dual-frame expansion established in the proof of
	Proposition~\ref{prop:densidad} converges in $\mathcal S$ and gives
	\[
	f=\sum_r a_r\phi_r.
	\]
	The local bounds of Lemma~\ref{lem:m-local} imply at most polynomial
	$L^2$ growth of $m$ on high-frequency annuli; hence
	$m\widehat f$ defines a tempered function and the usual Fourier multiplier
	$T_m f=\mathcal F^{-1}(m\widehat f)$ is well defined. For each frame index
	$w$, the compact Fourier support of $\widehat\phi_w$ and the convergence
	of the preceding expansion in $\mathcal S$ allow termwise pairing, yielding
	\[
	\begin{aligned}
	\langle T_m f,\phi_w\rangle
	&=
	\sum_r a_r\langle T_m\phi_r,\phi_w\rangle \\
	&=
	(\mathcal G^m a)_w.
	\end{aligned}
	\]
	Therefore
	\begin{equation}\label{eq:Tm-analysis-identity}
	C_\phi(T_m f)
	=
	\mathcal G^m D C_\phi f.
	\end{equation}
	The right-hand side belongs to $X_\beta$ by
	Lemma~\ref{lem:matrix-Tm}, so $T_m f\in M^{p,q}_{\alpha,\beta}$.
	Applying the retraction identity to \eqref{eq:Tm-analysis-identity} gives
	\[
	T_m f
	=
	R_\phi C_\phi(T_m f)
	=
	R_\phi\mathcal G^m D C_\phi f
	=
	\mathcal T_m f.
	\]
	Thus the coefficient-space operator agrees with the original Fourier
	multiplier on $\mathcal S$. We henceforth denote $\mathcal T_m$ again by
	$T_m$.

	If $p,q<\infty$, Proposition~\ref{prop:densidad} shows that
	$\mathcal S$ is dense in the source space, so this bounded extension is
	unique. For $p=\infty$ or $q=\infty$,
	\eqref{eq:Tm-coefficient-extension} provides the announced bounded
	endpoint realization.
\end{proof}

\begin{remark}[The Hilbertian case and relation with previous work]
	By Proposition~\ref{prop:L2},
	\[
	M^{2,2}_{\alpha,\beta}=H^\beta
	\]
	with equivalent norms. Consequently, the case $p=q=2$ of the theorem is Hilbertian in nature.
Indeed, Lemma~\ref{lem:m-local} together with the rescaled Sobolev
embedding, using $N_d>d/2$, yields
\[
|m(\xi)|
\lesssim
(1+|\xi|)^{-\mu}
\]
on the high-frequency region. Hence
\[
T_m:H^{\beta-\mu}\longrightarrow H^\beta
\]
follows directly from Plancherel's theorem. Thus the Hilbertian
estimate by itself is not the new content of the argument.

	Boundedness of Fourier multipliers on classical modulation,
	$\theta$--modulation, and more general decomposition spaces has an
	extensive literature; see, for instance,
	\cite{BenyiGrochenigOkoudjouRogers2007,
	ZhaoChenFanGuo2016,CleanthousGeorgiadisNielsen2018,
	GrochenigRzeszotnik2008}.
	Accordingly, we do not claim the general multiplier principle as new.
	For $0<\alpha\le1$, Proposition~\ref{prop:relation-alpha-mod} places the
	result within the existing $\theta$--modulation-space framework.
	The specific point here is the direct frame-matrix realization of the
	Beltr\'an--Bennett subdyadic geometry, including the shrinking-cell
	regime $\alpha>1$.
\end{remark}
\section{A subdyadic Gabor-frame characterization of the wavefront set}
\label{seccion:WFs:GaborSudDiadico}

We now use the subdyadic frame to give a phase-space characterization
of H\"ormander's classical local wavefront set. The construction below
should not be confused with the global Gabor wavefront set of
Rodino--Wahlberg \cite{RodinoWahlberg2014}, which is defined by decay in
conic subsets of the full phase space and coincides with H\"ormander's
global wavefront set. Likewise, the anisotropic Gabor wavefront sets of
\cite{RodinoWahlberg2023,Wahlberg2024} are global notions based on
power-type curves in phase space.

Here spatial localization is performed first by a compactly supported
cut-off, and the subdyadic frame is then used to measure conic
high-frequency regularity. Thus the resulting set is local in the
usual microlocal sense.

\subsection{Definition through subdyadic frame coefficients}

Fix the window $\phi$ and the subdyadic lattice
$\Omega_\alpha$ constructed in Proposition~\ref{prop:frame}. Thus
\[
\Omega_\alpha
=
\{(x_j^k,\xi_k):j\in\mathbb Z^d,\ k\in\mathbb N\},
\qquad
r_k:=|\xi_k|^{1-\alpha},
\]
with
\[
x_j^k=b\rho_\ast r_k^{-1}j,
\]
and
\[
\operatorname{supp}\widehat\phi\subset B(0,R_1),
\qquad
|\widehat\phi|\ge m_0
\quad\text{on }B(0,R_0),
\]
where
\[
A_{\mathrm{cov}}\rho_\ast<R_0,
\qquad
b\rho_\ast R_1<\pi.
\]
These are precisely the quantitative hypotheses used in the
Fourier-sampling proof of Proposition~\ref{prop:frame}.

For a compactly supported distribution $v\in\mathcal E'(\mathbb R^d)$
define the $k$-th block energy
\begin{equation}\label{eq:wf-block-energy}
\mathcal E_k^\phi(v)
:=
\left(
\sum_{j\in\mathbb Z^d}
|\langle v,\phi_{x_j^k,\xi_k}\rangle|^2
\right)^{1/2}.
\end{equation}
Since $\widehat v$ is smooth and of polynomial growth, the blockwise
Parseval argument used in Proposition~\ref{prop:frame} applies and gives
\begin{equation}\label{eq:wf-block-parseval}
\bigl(\mathcal E_k^\phi(v)\bigr)^2
=
C_{\rho_\ast}
\int_{\mathbb R^d}
|\widehat v(\eta)|^2
\left|
\widehat\phi
\left(
\frac{\eta-\xi_k}{r_k}
\right)
\right|^2
\,d\eta,
\end{equation}
where
\[
C_{\rho_\ast}
=
(2\pi)^{-d}(b\rho_\ast)^{-d}
\]
is independent of $k$.

\begin{definition}[Subdyadic Gabor microlocal regularity]
\label{def:reg-gabor}
Let $u\in\mathcal S'(\mathbb R^d)$ and
$(x_0,\xi_0)\in
\mathbb R^d\times(\mathbb R^d\setminus\{0\})$.
We say that $u$ is \emph{subdyadically Gabor-regular of type $\alpha$}
at $(x_0,\xi_0)$ if there exist
\[
\chi\in C_c^\infty(\mathbb R^d),
\qquad
\chi\equiv1
\quad\text{in a neighbourhood of }x_0,
\]
an open cone
$\Gamma\subset\mathbb R^d\setminus\{0\}$ containing $\xi_0$,
and $R>0$ such that, for every $N>0$, there exists $C_N>0$ with
\begin{equation}\label{eq:wf-block-decay}
\mathcal E_k^\phi(\chi u)
\le
C_N\langle\xi_k\rangle^{-N}
\end{equation}
whenever
\[
\xi_k\in\Gamma,
\qquad
|\xi_k|\ge R.
\]
\end{definition}

\begin{definition}[Subdyadic Gabor wavefront set]
\label{def:wf-gabor}
The \emph{subdyadic Gabor wavefront set} is
\[
WF_\alpha^{\mathrm G}(u)
=
\left\{
(x_0,\xi_0):
u\text{ is not regular at }(x_0,\xi_0)
\text{ in the sense of Definition~\ref{def:reg-gabor}}
\right\}.
\]
\end{definition}

The restriction $|\xi_0|^\alpha\ge1$ is deliberately absent from this
definition. Wavefront sets depend on directions, not on the magnitude
of the representative covector, while the estimates themselves concern
the limit $|\xi_k|\to\infty$. Since $\alpha>0$, the dispersive region
$|\xi_k|^\alpha\ge1$ is automatically reached at high frequency.

\subsection{Identification with the classical local wavefront set}

\begin{theorem}[Subdyadic characterization of the classical wavefront set]
\label{thm:wf-classical}
For every $\alpha>0$ and every
$u\in\mathcal S'(\mathbb R^d)$,
\[
WF_\alpha^{\mathrm G}(u)
=
WF(u),
\]
where $WF(u)$ denotes H\"ormander's classical $C^\infty$ wavefront set.
\end{theorem}

\begin{proof}
We prove both implications.

\medskip
\noindent\emph{Step 1: classical microlocal regularity implies
subdyadic coefficient decay.}
Suppose
\[
(x_0,\xi_0)\notin WF(u).
\]
Then there exist
$\chi\in C_c^\infty$ with $\chi\equiv1$ near $x_0$
and an open cone $\Gamma_0$ containing $\xi_0$ such that, for every
$M>0$,
\[
|\widehat{\chi u}(\eta)|
\le
C_M\langle\eta\rangle^{-M},
\qquad
\eta\in\Gamma_0.
\]
Choose a cone $\Gamma$ whose closure on the unit sphere is contained in
$\Gamma_0$.

If $\xi_k\in\Gamma$ and $|\xi_k|$ is sufficiently large, then
\[
\operatorname{supp}
\widehat{\phi_{x_j^k,\xi_k}}
\subset
B(\xi_k,R_1r_k)
\subset\Gamma_0.
\]
Indeed,
\[
\frac{r_k}{|\xi_k|}
=
|\xi_k|^{-\alpha}
\longrightarrow0.
\]
Moreover,
$\langle\eta\rangle\asymp\langle\xi_k\rangle$
on this ball. Hence \eqref{eq:wf-block-parseval} gives
\[
\bigl(\mathcal E_k^\phi(\chi u)\bigr)^2
\lesssim
r_k^d\langle\xi_k\rangle^{-2M}.
\]
Since
\[
r_k^d
=
|\xi_k|^{d(1-\alpha)},
\]
and $M$ is arbitrary, this implies
\[
\mathcal E_k^\phi(\chi u)
\le
C_N\langle\xi_k\rangle^{-N}
\]
for every $N>0$. Thus
$(x_0,\xi_0)\notin WF_\alpha^{\mathrm G}(u)$.

\medskip
\noindent\emph{Step 2: subdyadic coefficient decay implies classical
microlocal regularity.}
Conversely, suppose
\[
(x_0,\xi_0)\notin WF_\alpha^{\mathrm G}(u),
\]
and let $\chi$ and $\Gamma$ be as in
Definition~\ref{def:reg-gabor}. Put $v=\chi u$.

Since
\[
|\widehat\phi(\zeta)|
\ge m_0
\qquad
(|\zeta|\le R_0),
\]
formula \eqref{eq:wf-block-parseval} implies
\begin{equation}\label{eq:wf-inner-ball}
\int_{B(\xi_k,R_0r_k)}
|\widehat v(\eta)|^2\,d\eta
\lesssim
\bigl(\mathcal E_k^\phi(v)\bigr)^2.
\end{equation}

Choose an open cone $\Gamma'$ whose closure on the unit sphere is
contained in $\Gamma$. By the frequency covering of
Lemma~\ref{lem:covering} and
$A_{\mathrm{cov}}\rho_\ast<R_0$,
for every sufficiently large $\eta\in\Gamma'$ there exists $k$ such
that
\[
|\eta-\xi_k|
\le
A_{\mathrm{cov}}\rho_\ast r_k
<
R_0r_k.
\]
Since
$r_k/|\xi_k|=|\xi_k|^{-\alpha}\to0$,
the corresponding centre $\xi_k$ belongs to $\Gamma$ for all sufficiently
large $|\eta|$, and
\[
\langle\eta\rangle\asymp\langle\xi_k\rangle.
\]

Let $s\in\mathbb R$. Using the preceding covering and
\eqref{eq:wf-inner-ball},
\[
\int_{\Gamma',\,|\eta|\ge R'}
\langle\eta\rangle^{2s}
|\widehat v(\eta)|^2\,d\eta
\lesssim
\sum_{\substack{k:\,\xi_k\in\Gamma\\|\xi_k|\ge R''}}
\langle\xi_k\rangle^{2s}
\bigl(\mathcal E_k^\phi(v)\bigr)^2.
\]
The separation in Lemma~\ref{lem:covering} and a packing argument give
\[
\#\left\{
k:
2^j\le|\xi_k|<2^{j+1}
\right\}
\lesssim
2^{j\alpha d}.
\]
Using \eqref{eq:wf-block-decay} with $N$ arbitrarily large, the last
series is therefore bounded by
\[
C_N
\sum_{j\ge j_0}
2^{j\alpha d}\,
2^{2j(s-N)},
\]
which converges once
\[
N>s+\frac{\alpha d}{2}.
\]
Hence
\[
\langle\eta\rangle^s\widehat v(\eta)
\in
L^2(\Gamma')
\qquad
\text{for every }s\in\mathbb R.
\]
By the standard Sobolev characterization of H\"ormander's
$C^\infty$ wavefront set \cite{HormanderVol1},
\[
(x_0,\xi_0)\notin WF(u).
\]
This proves the equality.
\end{proof}

\subsection{Independence and localization}

\begin{proposition}[Robustness of the characterization]
\label{prop:wf-window}
Let $(\phi,\Omega_\alpha)$ and
$(\phi',\Omega'_\alpha)$ be two subdyadic frame pairs satisfying the
quantitative hypotheses used in Proposition~\ref{prop:frame}, namely
compact Fourier support, a uniform lower bound on an inner frequency
ball, and sufficiently fine painless sampling relative to the associated
subdyadic covering. Then both coefficient constructions define the same
microlocal set.
\end{proposition}

\begin{proof}
The proof of Theorem~\ref{thm:wf-classical} applies to either pair and
shows that both resulting sets coincide with H\"ormander's classical
wavefront set $WF(u)$.
\end{proof}

\begin{proposition}[Localization by cut-offs]
\label{prop:wf-cutoff}
Suppose that $\chi\equiv1$ on a neighbourhood $V$ of $x_0$ and that
the decay condition in Definition~\ref{def:reg-gabor} holds for
$\chi u$ in a cone around $\xi_0$. If
$\chi_1\in C_c^\infty(V)$, then the same microlocal regularity holds
for $\chi_1u$.
\end{proposition}

\begin{proof}
This is the usual localization property of the classical wavefront set,
combined with Theorem~\ref{thm:wf-classical}.
\end{proof}

\begin{proposition}\label{prop:wf-empty}
For $u\in\mathcal S'(\mathbb R^d)$,
\[
WF_\alpha^{\mathrm G}(u)=\varnothing
\quad\Longleftrightarrow\quad
u\in C^\infty(\mathbb R^d).
\]
In particular, if
\[
u\in
\bigcap_{\beta\in\mathbb R}
M^{2,2}_{\alpha,\beta},
\]
then $WF_\alpha^{\mathrm G}(u)=\varnothing$.
\end{proposition}

\begin{proof}
The first statement follows from
Theorem~\ref{thm:wf-classical}.
For the second, Proposition~\ref{prop:L2} gives
\[
M^{2,2}_{\alpha,\beta}=H^\beta
\]
with equivalent norms, and membership in all Sobolev orders implies
smoothness.
\end{proof}

\begin{remark}[No continuous inversion is used]
\label{rem:continuous-reconstruction}
The characterization above uses only the discrete subdyadic frame and
the blockwise Parseval identity \eqref{eq:wf-block-parseval}. In
particular, no classical fixed-window STFT inversion formula is invoked.
As observed earlier, such a formula does not automatically apply to the
continuous family whose dilation depends on $\xi$.
\end{remark}

\subsection{Pseudodifferential microlocality}

We use the standard classical class
$\Psi^0_{1,0}(\mathbb R^d)$. Thus a symbol
$a\in S^0_{1,0}$ satisfies
\[
|D_x^\gamma D_\xi^\delta a(x,\xi)|
\le
C_{\gamma,\delta}
\langle\xi\rangle^{-|\delta|}.
\]
We say that $a$ is elliptic at $(x_0,\xi_0)$ if there exist a
neighbourhood $U$ of $x_0$, an open cone $\Gamma$ containing $\xi_0$,
and constants $c,R>0$ such that
\[
|a(x,\xi)|
\ge c,
\qquad
x\in U,\quad
\xi\in\Gamma,\quad
|\xi|\ge R.
\]

\begin{theorem}[Microlocality and ellipticity]
\label{thm:wf-psido}
Let $A\in\Psi^0_{1,0}(\mathbb R^d)$ be properly supported. Then
\[
WF_\alpha^{\mathrm G}(Au)
\subset
WF_\alpha^{\mathrm G}(u).
\]
If $A$ is elliptic at $(x_0,\xi_0)$, then
\[
(x_0,\xi_0)\notin WF_\alpha^{\mathrm G}(u)
\quad\Longleftrightarrow\quad
(x_0,\xi_0)\notin WF_\alpha^{\mathrm G}(Au).
\]
\end{theorem}

\begin{proof}
By Theorem~\ref{thm:wf-classical},
\[
WF_\alpha^{\mathrm G}(u)=WF(u),
\qquad
WF_\alpha^{\mathrm G}(Au)=WF(Au).
\]
The assertions are therefore precisely the standard microlocality and
microellipticity properties of classical pseudodifferential operators
\cite{HormanderVol3}.
\end{proof}

\begin{remark}[Relation with Gabor and modulation-space wavefront sets]
The notation $WF_\alpha^{\mathrm G}$ records the analyzing frame and the
subdyadic resolution, not a new set of singular points:
\[
WF_\alpha^{\mathrm G}(u)=WF(u)
\qquad
\text{for every }\alpha>0.
\]
This is fundamentally different from the global Gabor wavefront set of
Rodino--Wahlberg \cite{RodinoWahlberg2014}, which is defined through
conic decay in the full phase space and coincides with H\"ormander's
global wavefront set. The later anisotropic Gabor wavefront sets
\cite{RodinoWahlberg2023,Wahlberg2024} are likewise global objects and
measure anisotropic behaviour along power-type curves in phase space.

Our construction is closer in spirit to the local wavefront-set theory
based on Fourier--Lebesgue and modulation spaces and on discrete Gabor
coefficients; see, for example,
\cite{JohanssonPilipovicTeofanovToft2012}.
For $0<\alpha<1$, the scale
$r_k\asymp|\xi_k|^{1-\alpha}$ is also consistent with variable-scale
wave-packet characterizations of the classical wavefront set such as
\cite{KatoKobayashiIto2017}. The point specific to the present paper is
the characterization through the explicit Beltr\'an--Bennett
subdyadic covering and its associated frame, including the shrinking
frequency-cell regime $\alpha>1$.
\end{remark}

\section*{Acknowledgments}
The author thanks the anonymous referee for a careful reading of the manuscript and for valuable comments that helped improve the presentation and strengthen several arguments.

	\bibliographystyle{abbrv}
	\bibliography{subdiadico_gabor_references}

@article{RodinoWahlberg2014,
	author = {L. Rodino and P. Wahlberg},
	title = {{The Gabor wave front set}},
	journal = {Monatshefte f{\"u}r Mathematik},
	volume  = {173},
	number  = {4},
	pages   = {625--655},
	year    = {2014},
	doi     = {10.1007/s00605-013-0592-0}
}

@article{RodinoWahlberg2023,
	author = {L. Rodino and P. Wahlberg},
	title = {{Anisotropic global microlocal analysis for tempered distributions}},
	journal = {Monatshefte f{\"u}r Mathematik},
	volume  = {202},
	number  = {2},
	pages   = {397--434},
	year    = {2023},
	doi     = {10.1007/s00605-022-01812-z}
}

@article{JohanssonPilipovicTeofanovToft2012,
	author = {K. Johansson and S. Pilipovi{\'c} and N. Teofanov and J. Toft},
	title = {{Gabor pairs, and a discrete approach to wave-front sets}},
	journal = {Monatshefte f{\"u}r Mathematik},
	volume  = {166},
	number  = {2},
	pages   = {181--199},
	year    = {2012},
	doi     = {10.1007/s00605-011-0288-2}
}

@article{KatoKobayashiIto2017,
	author = {K. Kato and M. Kobayashi and S. Ito},
	title = {{Remark on characterization of wave front set by wave packet transform}},
	journal = {Osaka Journal of Mathematics},
	volume  = {54},
	number  = {2},
	pages   = {209--228},
	year    = {2017},
	doi     = {10.18910/61907}
}

@article{DorflerMatusiak2014,
	author = {M. D{\"o}rfler and E. Matusiak},
	title = {{Nonstationary Gabor Frames---Existence and Construction}},
	journal = {International Journal of Wavelets, Multiresolution and Information Processing},
	volume  = {12},
	number  = {3},
	pages   = {1450032},
	year    = {2014},
	doi     = {10.1142/S0219691314500325}
}

@article{BorupNielsen2007,
	author = {L. Borup and M. Nielsen},
	title = {{Frame Decomposition of Decomposition Spaces}},
	journal = {Journal of Fourier Analysis and Applications},
	volume  = {13},
	number  = {1},
	pages   = {39--70},
	year    = {2007},
	doi     = {10.1007/s00041-006-6024-y}
}

@article{BenyiGrochenigOkoudjouRogers2007,
	author = {{\'A}. B{\'e}nyi and K. Gr{\"o}chenig and K. A. Okoudjou and L. G. Rogers},
	title = {{Unimodular Fourier multipliers for modulation spaces}},
	journal = {Journal of Functional Analysis},
	volume  = {246},
	number  = {2},
	pages   = {366--384},
	year    = {2007},
	doi     = {10.1016/j.jfa.2006.12.019}
}

@article{GrochenigRzeszotnik2008,
	author = {K. Gr{\"o}chenig and Z. Rzeszotnik},
	title = {{Banach algebras of pseudodifferential operators and their almost diagonalization}},
	journal = {Annales de l'Institut Fourier},
	volume  = {58},
	number  = {7},
	pages   = {2279--2314},
	year    = {2008},
	doi     = {10.5802/aif.2414}
}

@article{ZhaoChenFanGuo2016,
	author = {G. Zhao and J. Chen and D. Fan and W. Guo},
	title = {{Sharp estimates of unimodular multipliers on frequency decomposition spaces}},
	journal = {Nonlinear Analysis},
	volume  = {142},
	pages   = {26--47},
	year    = {2016},
	doi     = {10.1016/j.na.2016.04.003}
}

@article{CleanthousGeorgiadisNielsen2018,
	author = {G. Cleanthous and A. G. Georgiadis and M. Nielsen},
	title = {{Fourier Multipliers on Decomposition Spaces of Modulation and Triebel--Lizorkin Type}},
	journal = {Mediterranean Journal of Mathematics},
	volume  = {15},
	number  = {3},
	pages   = {122},
	year    = {2018},
	doi     = {10.1007/s00009-018-1171-3}
}

@phdthesis{Groebner1992,
	author = {P. Gr{\"o}bner},
	title = {{Banachr{\"a}ume glatter Funktionen und Zerlegungsmethoden}},
	school = {Universit{\"a}t Wien},
	year   = {1992}
}

@article{BorupNielsen2006,
	author = {L. Borup and M. Nielsen},
	title = {{Banach frames for multivariate {$\alpha$}-modulation spaces}},
	journal = {Journal of Mathematical Analysis and Applications},
	volume  = {321},
	number  = {2},
	pages   = {880--895},
	year    = {2006},
	doi     = {10.1016/j.jmaa.2005.08.091}
}

@article{Fornasier2007,
	author = {M. Fornasier},
	title = {{Banach frames for {$\alpha$}-modulation spaces}},
	journal = {Applied and Computational Harmonic Analysis},
	volume  = {22},
	number  = {2},
	pages   = {157--175},
	year    = {2007},
	doi     = {10.1016/j.acha.2006.05.008}
}

@article{Voigtlaender2022,
	author = {F. Voigtlaender},
	title = {{Structured, compactly supported Banach frame decompositions of decomposition spaces}},
	journal = {Dissertationes Mathematicae},
	volume  = {575},
	pages   = {1--179},
	year    = {2022},
	doi     = {10.4064/dm804-5-2021}
}

@article{FeichtingerGroebner1985,
	author = {H. G. Feichtinger and P. Gr{\"o}bner},
	title = {{Banach spaces of distributions defined by decomposition methods. {I}}},
	journal = {Mathematische Nachrichten},
	volume  = {123},
	number  = {1},
	pages   = {97--120},
	year    = {1985}
}

@article{FeichtingerGroebner1987,
	author = {H. G. Feichtinger},
	title = {{Banach spaces of distributions defined by decomposition methods. {II}}},
	journal = {Mathematische Nachrichten},
	volume  = {132},
	number  = {1},
	pages   = {207--237},
	year    = {1987}
}

@incollection{FeichtingerVoigtlaender2017,
	author = {H. G. Feichtinger and F. Voigtlaender},
	title = {{From {F}razier--{J}awerth characterizations of {B}esov spaces to wavelets and decomposition spaces}},
	booktitle = {Functional Analysis, Harmonic Analysis, and Image Processing: A Collection of Papers in Honor of Bj{\"o}rn Jawerth},
	series    = {Contemporary Mathematics},
	volume    = {693},
	publisher = {American Mathematical Society},
	address   = {Providence, RI},
	pages     = {185--216},
	year      = {2017}
}

@book{Lerner2010,
	author = {N. Lerner},
	title = {{Metrics on the phase space and non-selfadjoint pseudo-differential operators}},
	fseries = {Pseudo-Differential Operators. Theory and Applications},
	series = {Pseudo-Differ. Oper., Theory Appl.},
	issn = {2297-0355},
	volume = {3},
	isbn = {978-3-7643-8509-5},
	year = {2010},
	publisher = {Basel: Birkh{\"a}user},
	language = {English},
	zbMATH = {5622774},
	Zbl = {1186.47001}
}

@book{LinaresPonce2015,
	author = {F. Linares and G. Ponce},
	title = {{Introduction to Nonlinear Dispersive Equations}},
	series    = {Universitext},
	edition   = {2},
	publisher = {Springer},
	year      = {2015}
}

@article{KenigPonceVega1991,
	author = {C. E. Kenig and G. Ponce and L. Vega},
	title = {{Oscillatory integrals and regularity of dispersive equations}},
	journal = {Indiana University Mathematics Journal},
	volume  = {40},
	number  = {1},
	year    = {1991},
	pages   = {33--69}
}

@article{Wahlberg2024,
	author = {P. Wahlberg},
	title = {{Propagation of anisotropic {Gabor} wave front sets}},
	fjournal = {Proceedings of the Edinburgh Mathematical Society. Series II},
	journal = {Proc. Edinb. Math. Soc., II. Ser.},
	issn = {0013-0915},
	volume = {67},
	number = {3},
	pages = {674--698},
	year = {2024},
	language = {English},
	doi = {10.1017/S0013091524000269},
	zbMATH = {7930419},
	Zbl = {1558.46033}
}

@book{Heinonen2001,
	author = {J. Heinonen},
	title = {{Lectures on Analysis on Metric Spaces}},
	publisher = {Springer},
	series    = {Universitext},
	year      = {2001},
	address   = {New York},
	isbn      = {978-0387952321}
}

@article{Fefferman1972,
	author = {C. Fefferman and E. M. Stein},
	title = {{{{\(H^p\)}} spaces of several variables}},
	fjournal = {Acta Mathematica},
	journal = {Acta Math.},
	issn = {0001-5962},
	volume = {129},
	pages = {137--193},
	year = {1972},
	language = {English},
	doi = {10.1007/BF02392215},
	zbMATH = {3406634},
	Zbl = {0257.46078}
}

@article{Fefferman1970,
	author = {C. Fefferman},
	title = {{Inequalities for strongly singular convolution operators}},
	fjournal = {Acta Mathematica},
	journal = {Acta Math.},
	issn = {0001-5962},
	volume = {124},
	pages = {9--36},
	year = {1970},
	language = {English},
	doi = {10.1007/BF02394567},
	zbMATH = {3300878},
	Zbl = {0188.42601}
}

@article{Hirschman1959,
	author = {I. I. {Hirschman, Jr.}},
	title = {{On multiplier transformations}},
	fjournal = {Duke Mathematical Journal},
	journal = {Duke Math. J.},
	issn = {0012-7094},
	volume = {26},
	pages = {221--242},
	year = {1959},
	language = {English},
	doi = {10.1215/S0012-7094-59-02623-7},
	zbMATH = {3139142},
	Zbl = {0085.09201}
}

@book{Wainger1965,
	author = {S. Wainger},
	title = {{Special trigonometric series in $k$-dimensions}},
	fseries = {Memoirs of the American Mathematical Society},
	series = {Mem. Am. Math. Soc.},
	issn = {0065-9266},
	volume = {59},
	isbn = {978-0-8218-1259-4; 978-1-4704-0005-7},
	year = {1965},
	publisher = {Providence, RI: American Mathematical Society (AMS)},
	language = {English},
	doi = {10.1090/memo/0059},
	zbMATH = {3221645},
	Zbl = {0136.36601}
}

@article{BeltranBennett2017,
	author = {D. Beltr{\'a}n and J. Bennett},
	title = {{Subdyadic square functions and applications to weighted harmonic analysis}},
	journal = {Advances in Mathematics},
	volume  = {307},
	pages   = {72--99},
	year    = {2017},
	doi     = {10.1016/j.aim.2016.11.018}
}

@book{Stein1993Harmonic,
	author = {E. M. Stein},
	title = {{Harmonic Analysis: Real-Variable Methods, Orthogonality, and Oscillatory Integrals}},
	publisher = {Princeton University Press},
	address   = {Princeton, NJ},
	year      = {1993},
	series    = {Princeton Mathematical Series},
	volume    = {43}
}

@article{Miyachi1980Wave,
	author = {A. Miyachi},
	title = {{On some estimates for the wave equation in {$L^{p}$} and {$H^{p}$}}},
	journal = {Journal of the Faculty of Science, University of Tokyo. Section IA Mathematics},
	volume  = {27},
	number  = {2},
	pages   = {331--354},
	year    = {1980}
}

@article{Miyachi1981Singular,
	author = {A. Miyachi},
	title = {{On some singular {Fourier} multipliers for {$H^p(\mathbb{R}^n)$}}},
	journal = {Journal of the Faculty of Science, University of Tokyo. Section IA Mathematics},
	volume  = {28},
	number  = {2},
	pages   = {267--315},
	year    = {1981}
}

@book{HormanderVol1,
	author = {L. H{\"o}rmander},
	title = {{The Analysis of Linear Partial Differential Operators I: Distribution Theory and Fourier Analysis}},
	publisher = {Springer},
	address   = {Berlin},
	year      = {1983},
	series    = {Grundlehren der Mathematischen Wissenschaften},
	volume    = {256}
}

@book{HormanderVol3,
	author = {L. H{\"o}rmander},
	title = {{The Analysis of Linear Partial Differential Operators III: Pseudo-Differential Operators}},
	publisher = {Springer},
	address   = {Berlin},
	year      = {1985},
	series    = {Grundlehren der Mathematischen Wissenschaften},
	volume    = {274}
}

@article{FeichtingerGrochenig1989,
	author = {H. G. Feichtinger and K. Gr{\"o}chenig},
	title = {{Banach spaces related to integrable group representations and their atomic decompositions {I}}},
	journal = {Journal of Functional Analysis},
	volume  = {86},
	number  = {2},
	pages   = {307--340},
	year    = {1989}
}

@book{Grochenig2001,
	author = {K. Gr{\"o}chenig},
	title = {{Foundations of Time-Frequency Analysis}},
	publisher = {Birkh{\"a}user},
	address   = {Boston},
	year      = {2001},
	series    = {Applied and Numerical Harmonic Analysis}
}

@article{Gabor1946,
	author = {D. Gabor},
	title = {{Theory of communication}},
	journal = {Journal of the Institution of Electrical Engineers, Part I-III: Radio and Communication Engineering},
	volume  = {93},
	number  = {26},
	pages   = {429--457},
	year    = {1946}
}

@article{CorderoNicolaRodino2010,
	author = {E. Cordero and F. Nicola and L. Rodino},
	title = {{Time-frequency analysis of {Fourier} integral operators}},
	fjournal = {Communications on Pure and Applied Analysis},
	journal = {Commun. Pure Appl. Anal.},
	issn = {1534-0392},
	volume = {9},
	number = {1},
	pages = {1--21},
	year = {2010},
	language = {English},
	doi = {10.3934/cpaa.2010.9.1},
	zbMATH = {5716420},
	Zbl = {1196.42027}
}

@book{Stein1970SingularIntegrals,
	author = {E. M. Stein},
	title = {{Singular Integrals and Differentiability Properties of Functions}},
	publisher = {Princeton University Press},
	address   = {Princeton, NJ},
	year      = {1970},
	series    = {Princeton Mathematical Series},
	volume    = {30}
}

@article{RubioDeFrancia1985,
	author = {J. L. {Rubio de Francia}},
	title = {{A {Littlewood--Paley} inequality for arbitrary intervals}},
	journal = {Revista Matem{\'a}tica Iberoamericana},
	volume  = {1},
	number  = {2},
	pages   = {1--14},
	year    = {1985}
}
\end{document}